\documentclass[hidelinks,onefignum,onetabnum]{siamart250211}

\usepackage{lipsum}
\usepackage{amsfonts}
\usepackage{graphicx}
\usepackage{algorithmic}
\usepackage{xspace}
\usepackage{latexsym}
\usepackage{amsmath}
\usepackage{amssymb}
\usepackage{lineno}
\usepackage{comment}
\usepackage{enumerate} 
\usepackage{appendix}
\usepackage{algorithm}
\usepackage{tabularx}
\usepackage{listings}
\usepackage{multirow}
\usepackage{booktabs} 
\usepackage{ragged2e}
\usepackage{tikz,xcolor}
\usepackage{bm} 
\usepackage{amsopn}

\newsiamremark{remark}{Remark}
\newsiamremark{hypothesis}{Hypothesis}
\crefname{hypothesis}{Hypothesis}{Hypotheses}
\newsiamthm{claim}{Claim}
\newsiamremark{fact}{Fact}
\crefname{fact}{Fact}{Facts}

\headers{ST--SE SBP--SATs method for transient topology optimization}{S.~Nataj, M.~Appel, J.~Alexandersen}
\title{Space--time spectral element summation-by-parts method for heterogeneous transient heat conduction arising in topology optimization\thanks{
The work has been funded by Independent Research Foundation Denmark (DFF) through a Sapere Aude Research Leader grant (3123-00020B) for the COMFORT project (COmputational Morphogenesis FOR Time-dependent problems).}}

\author{%
  Sarah Nataj\thanks{Institute of Mechanical and Electrical Engineering, University of Southern Denmark, (\email{rena@sdu.dk}, \email{magap@sdu.dk}, \email{joal@sdu.dk}).}
  \and
  Magnus Appel\footnotemark[2]
  \and
  Joe Alexandersen\footnotemark[3]
}

\DeclareMathOperator{\diag}{diag}
\usetikzlibrary{shapes.geometric, arrows, positioning}

\tikzstyle{mainblock} = [rectangle, rounded corners, minimum width=3cm, minimum height=1cm,text centered, draw=black, fill=blue!30]
\tikzstyle{arrow} = [thick,->,>=stealth]

\definecolor{mygreen}{RGB}{44,85,17}
\definecolor{myblue}{RGB}{34,31,217}
\definecolor{mybrown}{RGB}{194,164,113}
\definecolor{myred}{RGB}{255,66,56}

\newcommand{\Rs}{\mathbb R}

\newcommand{\calA}{{\mathcal A}}

\newcommand{\calR}{{\mathcal R}}
\newcommand{\calL}{{\mathcal L}}

\renewcommand{\L}{{\mathcal{L}}}
\newcommand{\x}{{\boldsymbol{x}}}
\renewcommand{\u}{{\boldsymbol{u}}}
\newcommand{\f}{{\boldsymbol{f}}}
\newcommand{\g}{{\boldsymbol{g}}}
\newcommand{\y}{{\boldsymbol{y}}}
\newcommand{\z}{{\boldsymbol{z}}}

\newcommand{\vrho}{{\boldsymbol{\rho}}}

\newcommand{\zero}{\mbox{\boldmath $0$}}

\newcommand{\e}{\mbox{\boldmath $e$}}
\renewcommand{\b}{\mbox{\boldmath $b$}}

\renewcommand{\u}{\mbox{\boldmath $u$}}
\renewcommand{\v}{\mbox{\boldmath $v$}}

\newcommand{\q}{\mbox{\boldmath $q$}}
\newcommand{\h}{\mbox{\boldmath $h$}}

\newcommand{\Tr}{\ensuremath{^{\mathsf{T}}}} 
\newcommand{\mat}[1]{\ensuremath{\mathsf{#1}}} 

\newcommand{\DoneD}[1]{\ensuremath{\overline{\mat{D}}}_{#1}}
\newcommand{\QoneD}[1]{\ensuremath{\overline{\mat{Q}}}_{#1}}
\newcommand{\PoneD}[1]{\ensuremath{\overline{\mat{P}}}_{#1}}

\newcommand{\EoneD}[1]{\ensuremath{\overline{\mat{E}}}_{#1}}

\newcommand{\D}[1]{\ensuremath{\mat{D}}_{#1}}
\renewcommand{\S}[2]{\ensuremath{\mat{S}_{#1}^{#2}}}

\newcommand{\Q}[1]{\ensuremath{\mat{Q}}_{#1}}
\newcommand{\Pnorm}[1]{\ensuremath{\mat{P}}_{#1}}
\newcommand{\E}[1]{\ensuremath{\mat{E}}_{#1}}

\newcommand{\R}[1]{\ensuremath{\mat{R}}_{#1}}
\newcommand{\I}[1]{\ensuremath{\mat{I}}_{#1}}

\newcommand{\Jac}{\ensuremath{\mathsf{J}}}

\newtheorem{thm}{Theorem}

\newtheorem{lem}[thm]{Lemma}

\hypersetup{
	pdftitle={Space--time spectral element summation-by-parts method for heterogeneous transient heat conduction arising in topology optimization},
	pdfauthor={S.~Nataj, M.~Appel, J.~Alexandersen}
}

\begin{document}

\maketitle

\begin{abstract}
We develop a stable space--time spectral element method for transient heat conduction in heterogeneous multi-material domains, motivated by topology optimization. The method treats space and time simultaneously, using summation-by-parts (SBP) operators in both directions and simultaneous approximation terms (SATs) to impose boundary, initial, terminal, and material-interface conditions weakly, yielding a stable monolithic space–time scheme on heterogeneous domains. Stability is proven under specific conditions on the SAT parameters, scaled with the spatial mesh resolution and material properties. 
We compute design sensitivities using a discrete space–time adjoint scheme that is dual-consistent with the primal SBP–SAT scheme.  Dual consistency ensures that the discrete adjoint consistently approximates the continuous dual problem. From adjoint-consistency theory, under standard smoothness and stability assumptions, and for SBP-compatible functional discretization, dual consistency is expected to yield superconvergent functional estimates. 
This means that, for a fixed design, the error in the computed objective functional converges at a higher rate than the error in the underlying state approximation. In the topology optimization setting, this is relevant because the design update is driven by repeated evaluations of the objective functional and its adjoint-based sensitivity. Improved functional accuracy can therefore contribute to more reliable design updates. Numerical experiments demonstrate spectral convergence of the forward and adjoint solution errors and of the error in the functional output for smooth manufactured solutions in heterogeneous domains. Also under mesh refinement, the observed errors indicate higher-order convergence of the optimal objective value compared with the optimized design. We validate the resulting optimal design by comparison with an independently computed reference optimal design and report time-to-solution and cost-of-accuracy curves, comparing against low-order time-marching and all-at-once solvers for the forward and adjoint systems. The proposed scheme attains high accuracy with fewer space–time degrees of freedom and remains stable, reducing time-to-solution and memory compared with an alternative all-at-once solver. This makes it a future candidate for large-scale topology optimization of time-dependent thermal systems.



\end{abstract}


\section{Introduction}\label{sec:intro}

\subsection{Motivation}\label{sec:intro:motv}
Topology optimization is a computational design methodology for finding an optimal material distribution within a given domain that minimizes a given objective functional subject to governing physics and design constraints~\cite{Bendsoe1988,Sigmund2004}. Beyond classical structural applications in solid mechanics, modern uses span electromagnetic, thermal, acoustic, and fluid dynamics problems; see the survey articles~\cite{Deaton2014,Alex2020}. Here we focus on transient heat conduction and optimize material layouts to steer the temperature field over time. Although transient heat-conduction topology optimization has been studied in the literature~\cite{Zhuang2013,Zhuang2014,Dbouk2017,Long2018,Wu2019,Zeng2020,Li2023}, the use of high-order discretizations with stability and adjoint-consistency properties remains limited.

We adopt an adjoint-based sensitivity framework for the partial differential equation (PDE) constrained optimization problem.  Design sensitivities are computed via a discrete adjoint and require one forward and one backward solve. An overview of density-based and alternative topology optimization approaches is given  in \cite{Sigmund2013}. Unlike a topology optimization problem for steady-state PDEs, where only one state needs to be computed and stored for adjoint sensitivity analysis, an optimization problem subject to time-dependent PDEs must consider an entire time history \cite{Michaleris1994,Dahl2008}.
Therefore, the adjoint method for unsteady problems typically requires storing the state at all time steps. This, plus the necessity to perform many optimization iterations, makes time-dependent topology optimization extremely computationally demanding.
Another fundamental challenge is the sequential nature of time integration. Conventional approaches march forward in time step by step for the state, then backward in time for the adjoint. This limits parallelism, even if the spatial domain is parallelized among processors, as computations  must   be performed sequentially in time. 

\subsection{Related work}\label{sec:intro:back}
Space–time methods provide a powerful approach to address these challenges. In a space–time formulation, the time dimension is treated similarly to a spatial dimension, enabling a fully simultaneous solution over all time steps. This approach eliminates the sequential time-stepping bottleneck by solving all time steps in a single unified system, preferably by parallelizing across time. The trade-off is that the associated linear system is significantly larger and often ill-conditioned. Efficient parallel solvers and preconditioners are required to handle the resulting very large linear systems. 
Recent literature has begun to explore strategies in topology optimization to reduce runtime \cite{Theulings2024,Appel2024} and memory usage \cite{Theulings2024,Margetis2023,Yaji2018} for transient problems.

Much of the work in this direction, however, employs low-order spatial and temporal schemes. 
Theulings et al. \cite{Theulings2024} use a finite-volume discretization on Cartesian grids with an explicit forward Euler scheme in time for solving a heat conduction problem and implicit backward Euler scheme for solving a fluid dynamics problem. Appel and Alexandersen’s one-shot Parareal study \cite{Appel2024} adopts a Galerkin finite-element discretization with bilinear spatial shape functions and first-order backward Euler scheme in time.  Their space–time multigrid paper \cite{Appel2025} targets an all-at-once BE–FE formulation (linear finite element in space, backward Euler in time) on uniform space–time meshes. In \cite{Alex2025} they use a continuous Galerkin space–time finite element method with linear elements in both space and time, with an added artificial diffusion in the time direction for stability, achieving a first-order accuracy in time. Overall, these studies emphasize solver efficiency and parallelism while retaining low-order discretizations. Low-order spatial elements are common in topology optimization, since the design freedom is coupled to the mesh resolution (usually one design variable per element) and thus requires many spatial elements to represent complicated structure. Low-order elements allow to keep total computational cost relatively low when needing fine meshes for high design freedom. However, this does not need to translate to temporal schemes, where usually the design will remain constant over time. Thus, we explore high-order temporal discretization herein.

High-order discretization of time-dependent PDEs in both space and time offers improved solution accuracy and convergence rates with fewer space--time degrees of freedom (DoFs) in comparison to standard low-order finite elements and finite difference methods.  Adjoint methods compatible with high-order discretizations for unsteady flow optimization have been established in the literature~\cite{Yamaleev2010,Zahr2016}, employing spectral element, Fourier spectral, or discontinuous Galerkin-type discretizations~\cite{Saglietti2018,Nobis2022,Nobis2023}. This provides the methodological basis for the present work, where the high-order space--time adjoint framework is extended to the  multi-material domains setting arising in topology optimization.

Schemes built from high-order differentiation operators that satisfy the summation-by-parts (SBP) property~\cite{Svard2014,DelRey2014} provide a principled 
path to stability in the fully discrete setting. Coupled with simultaneous approximation terms (SAT)~\cite{Carpenter1999,Mattsson2003},
they weakly enforce boundary and interface conditions while replicating integration-by-parts at the discrete level. The resulting SBP–SAT formulations admit discrete energy estimates that mirror their continuous problems and thereby prevent the spurious growth of numerical energy. The SBP-SAT approach has previously been used successfully to construct provably stable high-order discretizations of the conjugate heat transfer (CHT) problem using a semi-discrete approach~\cite{Lindstrom2010,Nordstrom2013,Lundquist2018}. The present work builds on the interface-coupling approach and stability analysis developed in~\cite{Nataj2026CHT}, extending the SBP--SAT framework from semi-discrete conjugate heat-transfer formulations to a fully space--time discretization of the transient heat equation in multi-material domains.

Building on SBP--in--time developments~\cite{Nordstrom2013time,Nordstrom2016,Lucas2019} 
and dual-consistent SBP schemes~\cite{Berg2012,Nikkar2019}, we construct a space–time spectral element SBP–SAT method for both primal and adjoint problems. Although the present work uses a spectral element (SE) discretization, the construction is expressed in terms of SBP operators and can therefore be adapted to other discretizations, such as finite-difference (FD), finite-volume (FV), discontinuous-Galerkin (DG), or continuous-Galerkin (CG) methods, provided that the required SBP structure is available.  The resulting framework treats high-order space--time discretization, weak interface coupling, adjoint consistency, and SBP-compatible output functionals in a unified way for the multi-material domain arising in topology optimization. Due to the SBP property of the differentiation operator together with compatibly constructed space/time boundary and interface SATs, the algebraic adjoint of the discrete state operator coincides with a consistent discretization of the continuous adjoint (discrete dual consistency). 
Dual (adjoint) consistency is a key property of a discretization and is typically required to obtain optimal-order convergence for given functionals \cite{Arnold2001,Hartmann2007}. In particular, a linear functional is called superconvergent if the error in its discrete approximation converges at a higher order than the underlying solution error \cite{Niles2000,Giles2002}. It is known that dual-consistent SBP finite-difference discretizations, when paired with an energy-stable discretization, yield superconvergent linear functionals under standard regularity conditions
~\cite{Hicken2011,Berg2013,Worku2022}. 
Motivated by these results, the present space--time spectral element SBP--SAT scheme is constructed to combine energy stability, dual consistency, and an SBP-compatible discretization of the output functional. Under the usual smoothness and stability assumptions, this combination is expected to improve the accuracy of functional estimates and the numerical results are consistent with this expectation.

\subsection{Contributions and scope}\label{sec:intro:cont}
The main contributions of this article are as follows.
\begin{enumerate}
    \item We develop a stable space--time spectral element SBP--SAT formulation for transient heat conduction in heterogeneous multi-material domains, as they arise in density-based topology optimization. In this formulation, the spatial and temporal directions are discretized simultaneously in a single monolithic space--time system, while material interfaces are treated with SAT penalty terms that weakly enforce continuity of temperature and normal conductive heat flux across neighbouring subdomains.

    \item We prove an energy estimate for the resulting monolithic space--time primal scheme under suitable choices of the SAT parameters, with scalings depending on the spatial mesh resolution and material properties. 

    \item We derive a discrete space--time adjoint that is dual-consistent with the primal SBP--SAT discretization, so that the discrete adjoint consistently approximates the continuous dual problem. Together with an SBP-compatible discretization of the output functional, this provides the setting in which improved functional accuracy is expected under standard smoothness and stability assumptions. 

    \item We apply the proposed construction in a density-based topology optimization loop for transient heat conduction and compare the proposed space--time spectral element formulation with a low-order all-at-once backward Euler finite-element solver.  
\end{enumerate}
 
To keep the notation and energy estimate transparent, the space-time stability analysis presented here is restricted to one spatial dimension. The tensor-product SBP--SAT construction, however, provides a natural route toward multidimensional multi-material problems with element-aligned interfaces. In higher dimensions, the same stability mechanism is expected to apply, but the proof involves additional geometric and notational complexity \cite{Nataj2026CHT}. The present paper therefore focuses on the one-dimensional space-time stability analysis, the dual-consistent discrete adjoint construction, and numerical evidence for improved functional accuracy. A formal functional-superconvergence proof for the space--time spectral-element setting is beyond the scope of the present analysis.

The numerical experiments use a one-dimensional setting as a controlled benchmark to validate the resulting optimized designs against an independently computed reference solution. They demonstrate spectral convergence of the forward and adjoint solution errors and of the functional-output error for smooth manufactured solutions in heterogeneous media. In the topology optimization setting, mesh-refinement studies further indicate higher-order accuracy in the optimal objective value compared with the optimized design. The results also show that the proposed approach reaches a prescribed accuracy with fewer temporal degrees of freedom and reduced time-to-solution. Extensions to higher-dimensional and large-scale optimization problems are a natural next step, but beyond the goal of the present work; these require additional treatment of filtering and projection, parallel implementation, and scalable solver strategies.

The paper is organized as follows. Section~\ref{sec:method} states the model and optimization problem, introduces the space-time spectral element SBP–SAT discretization,  proves the conditional stability of the space-time SBP-SAT scheme for heterogeneous heat conduction problem arising in topology optimization, and derives the discrete adjoint and sensitivities  formulae. Section~\ref{sec:results} presents numerical experiments that validate the optimal designs by comparison with independently computed reference designs and  assess accuracy and efficiency of the proposed scheme. Section~\ref{sec:sum} gives a conclusion.

\section{Methodology}\label{sec:method}
In this section, we formulate the high-order space--time  framework for a transient heat-conduction topology optimization problem.
For simplicity, this work considers heterogeneous transient heat conduction in one spatial dimension. Let $a<b$ be given and let $T>0$ denote the terminal time, and define the space--time domain $\Omega := (a,b)\times(0,T)$. Let $u(x,t)$ denote the temperature field and $f(x,t)$ the external heat load per unit volume.
The optimization problem is then formulated as a minimization of the squared temperature integrated over time and space to
reduce the temperature throughout the interval, subject to the transient heat conduction equation and a material volume bound \(V^*\). In contrast to the thermal compliance objective, which emphasizes regions where heat is applied, and therefore it is sensitive to the magnitude and localization of the source term, the chosen objective functional penalizes elevated temperature throughout the space--time domain. This type of temperature-control objective is relevant in transient thermal topology optimization, for example when reducing peak temperatures over a prescribed time interval~\cite{Wu2019}.

The presence of conductive or insulating material is parameterized by material distribution  function, $\rho(x)$, which is referred to as the ``design field'' and as ``design variables'' after discretization. Inside the insulator, $\rho(x)=0$, and inside the conductor, $\rho(x)=1$. 
In density-based topology optimization, $\rho$ is relaxed to vary continuously between $0$ and $1$, and penalized material interpolation is commonly employed to drive the solution toward the binary discrete values~\cite{Sigmund2004}. This is done by assigning intermediate densities a reduced effective material property, which makes them less favourable in the optimization process \cite{Sigmund2004}. Penalization alone, however, does not impose a binary constraint on the design variable.

Our continuous topology optimization problem is formulated as 
{\setlength{\abovedisplayskip}{6pt}
 \setlength{\belowdisplayskip}{6pt}
\begin{align}\label{eq00}
\begin{aligned}
	\underset{  u  , \,\rho}{\textrm{min.}}&\quad J(  u  ,\rho)
	= \int_{0}^{T }\!\!\int_{a}^{b}   u  (x,t)^2 \,\mathrm{d}x\,\mathrm{d}t,\\
	\text{s.t.\,:} &\quad \dfrac{\partial   u  }{\partial t}(x,t)
	= \dfrac{\partial}{\partial x}\!\Bigl(\kappa(\rho(x))\,
	\dfrac{\partial   u  }{\partial x}(x,t)\Bigr)
	+ f(x,t),
	\quad &(x,t)\in\Omega,\\
	&\quad   u  (a,t) 
	= h(t),\quad   u  (b,t) 
	= g(t),
	& t\in[0,T ],\\
	&\quad   u  (x,0) = q(x),
	& x\in[a,b],\\[0.75ex]
	&\int_{a}^{b}\rho(x)\,\mathrm{d}x \;\le\; V^*, \quad 0 \le \rho(x) \le 1,
\end{aligned}
\end{align}}
The thermal diffusivity \(\kappa\) depends on the material distribution $\rho(x)$, and is given by
 \begin{equation}\label{eq0}
 \kappa(\rho(x)) = \kappa_{\min} + (\kappa_{\max} - \kappa_{\min})\,\rho(x)^{p},\qquad 0 \le \rho(x) \le 1,
 \end{equation}
where \(\kappa_{\min}\), \(\kappa_{\max}\), and \(p\ge1\) are known positive constants. We consider a model problem in which the thermal diffusivity \(\kappa\) is defined as \(k/(\rho_m c_p)\), where \(k\) is the thermal conductivity, \(\rho_m\) is the mass density, and \(c_p\) is the specific heat capacity at constant pressure. 

We adopt a discretize--then--optimize strategy; the PDE and the objective are first discretized in space-time, where we employ SBP differentiation operators to discretize spatial and temporal derivatives and apply SATs to impose initial and interface conditions in a stable monolithic manner~\cite{Giles2000,Carpenter2010}.

In the remainder of this section, we detail these steps. Section~\ref{sec:sbp} reviews the basics of SBP operators; Section~\ref{sec:res} presents the assembly of the discrete forward residual using a  monolithic space-time spectral element SBP-SAT and Section~\ref{sec:stab} provides a proof of conditional stability of the proposed scheme. Section~\ref{sec:TopOpt} formulates the discrete topology optimization problem and describes the solution strategy; 
Section~\ref{sec:adj} derives the dual-consistent discrete adjoint scheme; Section~\ref{sec:sen} presents the calculation of the corresponding design sensitivities.

\subsection{SBP operators}\label{sec:sbp}

We employ SBP differentiation operators to discretize the spatial and temporal derivatives. Applying these SBP operators allows the use of high-order schemes and permits flexibility in the discretization type in each subdomain (e.g., DG, CG, FD, FV) \cite{Svard2014,DelRey2014}. 

Let $\{x_i\}_{i=0}^{N_x}$ denote the grid points on a given interval $[a, b]$. The one-dimensional SBP differentiation  operator of degree $r$, \(\DoneD{ x }\in \Rs^{(N_x+1)\times (N_x+1)}\), approximates the spatial derivative  \(\partial/ \partial  x ,\) and satisfies the following properties: it differentiates exactly  monomials up to degree $r$,
\[\DoneD{x} \, \x^s=s\,\x^{s-1}, \qquad s=0,1,\ldots, r,\qquad \x=[x_0,\ldots,x_{N_x}]^{\Tr},\]	
and it admits the factorization $\DoneD{ x }=\PoneD{ x }^{\,-1}\QoneD{ x }$,  where $\PoneD{ x }$ is symmetric positive definite and $\QoneD{ x }+\QoneD{ x }\Tr=\EoneD{ x }=\diag(-1,0,\ldots,0,1)$.

Next, we seek to construct space–time SBP differentiation operators 
together with a compatible quadrature rule 
for the space–time integrals, in a way that preserves discrete integration-by-parts identity. To extend the one-dimensional SBP operator to space–time, we first construct a one-dimensional temporal SBP operator 
$\DoneD{t} \approx \partial / \partial t$ of the form 
$\DoneD{t} = \PoneD{t}^{-1}\QoneD{t}$ on collocation nodes 
$\{t_j\}_{j=0}^{N_t} \subset [0,T]$. 
The space–time differentiation operators are then defined as tensor products \cite{DelRey2018}.
\begin{align*}
	\begin{aligned}	\Pnorm{}:=\PoneD{t}\otimes\PoneD{x},\qquad\Q{x}:=\PoneD{t}\otimes\QoneD{x},\qquad\quad\Q{t}:=\QoneD{t}\otimes\PoneD{x},
	\end{aligned}
\end{align*}
and then define
\begin{align*}
	\begin{aligned}
		 \D{x}:=\Pnorm{}^{-1}\Q{x}= \I{t}\otimes\DoneD{x },\qquad\D{t}:=\Pnorm{}^{-1}\Q{t}= \DoneD{t }\otimes \I{x},\qquad
	\end{aligned}
\end{align*}
where $\I{x}$ and $\I{t}$ are identity matrices of size $(N_x+1) \times (N_x+1)$ and $(N_t+1)\times (N_t+1)$, respectively. The symmetric positive definite weight (quadrature) matrix $\Pnorm{}$ defines an inner product via $\langle \x, \y\rangle_{\Pnorm{}}:=\x^T{\Pnorm{}}\y$ for given vectors $\x\in\Rs^{N}$ and $\y\in\Rs^{N}$, where $N=(N_x+1)(N_t+1)$. It also induces the discrete norm defined by ${\|\x\|_{\Pnorm{}}}^2=\x^T {\Pnorm{}}\x$. For vectors \(\u\) and \(\v\) containing the values of continuous functions  \( u (x,t)\) and \( v (x,t)\) at the grid points, respectively, one has 
\begin{align*} 
	\begin{aligned}			
    \u\Tr\Pnorm{}\v\approx\int_{\Omega} u  v \,d \Omega,
	\qquad\quad
	\u\Tr\Q{ x }\v\approx\int_{ \Omega} u \dfrac{\partial  v }{\partial x }\,d \Omega,
	\end{aligned}
\end{align*}
where $\Omega$ is the space-time domain.

We proceed to define restriction operators to extract the boundary values on the west (left), east (right), south (initial) and north (terminal) faces of each element:
\[\R{w} = \I{t} \otimes\e_w^\mathsf{T}, \qquad
\R{e} = \I{t} \otimes\e_e^\mathsf{T} , \qquad
\R{s} = \e_s^\mathsf{T}\otimes\I{x},  \qquad
\R{n} = \e_n^\mathsf{T}\otimes\I{x},
\]
where 
\begin{align*}
    \begin{aligned}
       &\e_w = [1,0,\dots,0]^\mathsf{T}\in \Rs^{N_x+1}, &\quad& \e_e = [0,\dots,0,1]^\mathsf{T}\in \Rs^{N_x+1},\\
       &\e_s = [1,0,\dots,0]^\mathsf{T}\in \Rs^{N_t+1}, &\quad& \e_n = [0,\dots,0, 1]^\mathsf{T}\in \Rs^{N_t+1}. 
    \end{aligned}
\end{align*}
The surface operators in the $x$ and $t$ directions are given as
{\setlength{\abovedisplayskip}{7pt}
 \setlength{\belowdisplayskip}{7pt}
\begin{align*} 
\begin{aligned}
\E{x}:=\R{e}\Tr\PoneD{t}\R{e}-\R{w}\Tr\PoneD{t}\R{w}, \qquad\E{t}:=\R{n}\Tr\PoneD{x}\R{n}-\R{s}\Tr\PoneD{x}\R{s} .
\end{aligned}
\end{align*}}
Notice that 
\begin{align*} 
	\begin{aligned}			
    \u\Tr\E{x}\v\approx\int_{\Gamma_x} u\, v\,  n_{ x }\,d{ \Gamma},\qquad
    \u\Tr\E{t}\v\approx\int_{\Gamma_t} u\, v\,  n_{ t }\,d{ \Gamma},
	\end{aligned}
\end{align*}
where $\Gamma_x$ and $\Gamma_t$ denote the spatial and temporal parts of the boundary of the space–time domain $\Omega$, respectively.

With these space–time SBP operators in place, we can now state the discrete integration-by-parts property that is central to the SBP framework. Consider integration-by-parts in continuous form for spatial dimension $x$,  
\begin{align*} 
	\begin{aligned}
		\int_{\Omega} u \,\dfrac{\partial  v }{\partial x }\, d\Omega+	\int_{ \Omega} v \,\dfrac{\partial  u }{\partial x }\,d \Omega=\int_{\Gamma_x} u\,  v\,  n_{ x }\,d{\Gamma}.
	\end{aligned}
\end{align*}
In discrete form we have the analogous property
\begin{align*} 
		\begin{aligned}
				\langle \u,\D{ x }\v\rangle_{\Pnorm{}}+\langle \v, \D{ x }\u\rangle_{\Pnorm{}}
				&=\u\Tr\Pnorm{}\D{ x }\v+\v\Tr\Pnorm{}(\D{ x }\u)\\
				&=\u\Tr\Pnorm{}(\Pnorm{}^{-1}\Q{ x })\v+\u\Tr (\Q{ x }\Tr\Pnorm{}^{-1})\Pnorm{}\v\\
				&=\u\Tr (\Q{ x }+\Q{ x }\Tr)\v=\u\Tr\E{ x }\v.
		\end{aligned}
\end{align*}

By replicating integration--by--parts at the discrete level, SBP-SAT schemes enforce boundary conditions weakly,  systematically creating  energy-stable numerical methods in conjugate heat transfer problems~\cite{DelRey2014,Svard2014}.

Spectral element methods based on Legendre–Gauss–Lobatto (LGL) nodes combine the geometric flexibility of finite element methods with the exponential convergence of spectral methods for smooth solutions. Their built-in SBP differentiation operators ensure provable energy stability, making them particularly well suited for CHT problems \cite{Gassner2013}. 

A space–time spectral element discretization for multi–domain CHT problems attains high-order in both space and time. Within the SBP–SAT framework, initial/terminal, boundary, and interface conditions are imposed weakly yet stably via SAT penalties, yielding a consistent scheme admitting a global energy estimate that mirrors the continuous energy estimate.

\subsection{Space-time spectral element discretization of heterogeneous transient heat conduction}\label{sec:res}
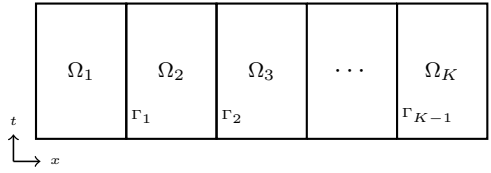
\begin{figure}[htp]
    \begin{Center}
    	\begin{tikzpicture}[scale=0.6]
    		\def\W{2}
    		\def\H{1.5}
    		\def\k{5}
    		\pgfmathsetmacro\lasti{\k-1}
    		\pgfmathsetmacro\lastx{\lasti*\W}
    		\pgfmathsetmacro\ellx{(3+\lasti)/2*\W}
    		\foreach \i in {1,2,3} {
    			\pgfmathsetmacro\x{(\i-1)*\W}
    			\draw[thick] (\x,-\H) rectangle ++(\W,2*\H);
    			\node at (\x+0.5*\W,0) {{\footnotesize$\Omega_{\i}$}};
    		}
    		\draw[thick] (3*\W,\H) -- (\lastx,\H);
    		\draw[thick] (3*\W,-\H) -- (\lastx,-\H);
    		\node at (\ellx,0) {$\cdots$};
    		\draw[thick] (\lastx,-\H) rectangle ++(\W,2*\H);
    		\node at (\lastx+0.5*\W,0) {{\footnotesize$\Omega_{K}$}};
    		\node[anchor=west] at (1*\W-0.1,-1)   {\tiny$\Gamma_{1}$};
    		\node[anchor=west] at (2*\W-0.1,-1)   {\tiny$\Gamma_{2}$};
    		\node[anchor=west] at (\lastx-0.1,-1) {\tiny$\Gamma_{K-1}$};
    		\coordinate (orig) at (-0.5,-\H-0.5);
    		\draw[->,line width=0.6pt] (orig) -- ++(0.6,0)  node[right] {\tiny$x$};
    		\draw[->,line width=0.6pt] (orig) -- ++(0,0.6)  node[above] {\tiny$t$};
    	\end{tikzpicture}
    \end{Center}	
    \caption{Illustration of the partitioning of the space-time domain, $\Omega$, into non-overlapping elements, $\Omega_k$, with interfaces, $\Gamma_k$.}
    \label{fig:domainElements}
\end{figure}

In the following, we construct a space–time SBP–SAT spectral element scheme for solving the primal problem. Suppose the given space-time domain $\Omega$ is partitioned into $K$ non-overlapping elements as illustrated in Figure \ref{fig:domainElements},
\[\bar\Omega \;=\;\bigcup_{k=1}^K \bar\Omega_k\;=\;\bigcup_{k=1}^K [a_k, b_k]\times [0,T].\]

Let $u:\bar\Omega\to\mathbb{R}$ denote the (continuous) solution of the heat equation and define $u^{(k)} := u|_{\Omega_k}$ as its restriction to element $\Omega_k$. Assume the elements are coupled at $\Gamma_k=\partial \Omega_{k}\cap \partial \Omega_{k+1}$ for $ 1 \le k \le K-1$, where continuity of the temperature  and continuity of the normal heat flux are enforced:
\begin{align}\label{eqint}
	\begin{aligned}
		  u  ^{(k)} &=   u  ^{(k+1)}, \quad\quad\qquad\quad \;\;\text{on } \Gamma_k, \\
		\kappa^{(k)}\,   u  _x^{(k)} &= \kappa^{(k+1)}\,   u  _x^{(k+1)}, \qquad\quad \text{on } \Gamma_k. \\
	\end{aligned}
\end{align}
In the standard setting of density-based topology optimization, the design variable represents a material distribution in the physical domain. For the transient problems, this material layout is fixed during one transient simulation, while the temperature field evolves in space and time \cite{Sigmund2004}. We therefore assume the thermal diffusivity in each element is constant, $\kappa^{(k)}=\kappa(\rho^{(k)})$ where $\rho^{(k)}=\rho(\bar{x})$ for some $\bar{x} \in \Omega_k$ for $ 1 \le k \le K$.
This ensures both energy conservation and the correct coupling at the interfaces between elements.

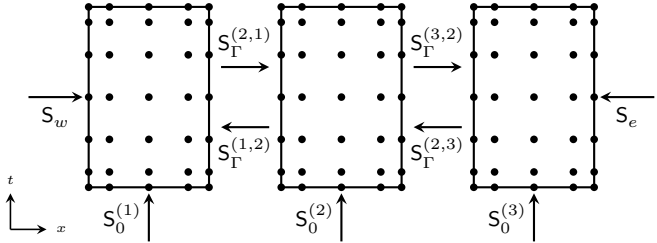
\begin{figure}[htp!]
\begin{Center}
\def\xLGL{0,0.345,1.000,1.655,2.000}
\def\yLGL{0,0.255,0.797,1.500,2.203,2.745,3.000}
\begin{tikzpicture}[scale=0.8, >=stealth]

\newcommand{\Subdomain}[2]{%
  \begin{scope}[shift={(#1,0)}]
    \draw[thick] (0,0) rectangle (2,3);
    \foreach \x in \xLGL {
      \foreach \y in \yLGL {
        \fill (\x,\y) circle (1.8pt);
      }
    }
    \draw[thick,->] (1,-0.9) -- (1,-0.1)
      node[below, xshift=-10pt] {{\footnotesize$\S{0}{(#2)}$}};
  \end{scope}%
}

\Subdomain{0}{1}
\Subdomain{3.2}{2}
\Subdomain{6.4}{3}

\draw[thick,->] (-1,1.5) -- (-0.1,1.5)
  node[below,pos=0.5] {{\footnotesize$\S{w}{}$}};

\draw[thick,<-] (8.5,1.5) -- (9.4,1.5)
  node[below,pos=0.5] {{\footnotesize$\S{e}{}$}};

\draw[thick,->] (2.2,2.0) -- (3.0,2.0)
  node[above,pos=0.5] {{\footnotesize$\S{\Gamma}{(2,1)}$}};
\draw[thick,->] (3.0,1.0) -- (2.2,1.0)
  node[below,pos=0.5] {{\footnotesize$\S{ \Gamma}{(1,2)}$}};

\draw[thick,->] (5.4,2.0) -- (6.2,2.0)
  node[above,pos=0.5] {{\footnotesize$\S{\Gamma}{(3,2)}$}};
\draw[thick,->] (6.2,1.0) -- (5.4,1.0)
  node[below,pos=0.5] {{\footnotesize$\S{\Gamma}{(2,3)}$}};

\coordinate (orig) at (-1.3,-0.7);
\draw[->,line width=0.6pt] (orig) -- ++(0.6,0)  node[right] {\tiny $x$};
\draw[->,line width=0.6pt] (orig) -- ++(0,0.6)  node[above] {\tiny $t$};

\end{tikzpicture}
\end{Center}
\caption{Illustration of three subdomains discretized using LGL collocation nodes and coupled using boundary, initial, and interface SATs.}
\label{fig:domaincoupling}
\end{figure}

To construct the SBP-SAT scheme, we assume $\DoneD{x}$ and  $\DoneD{t}$ are pseudospectral differentiation matrices based on LGL collocation nodes at element $k$,  
\[\{x_i^{(k)}\}_{i=0}^{N_x} \subset [a_k,b_k],\qquad \{t_j\}_{j=0}^{N_t} \subset [0,T].\]  
Define
\begin{align*}
\begin{aligned}
\u^{(k)} &:=\big[  u^{(k)}(x_0,t_0),
u^{(k)}(x_1,t_0),\ldots,
u ^{(k)} (x_{N_x},t_0),\\
&\qquad\qquad \dots,  u^{(k)}(x_0,t_{N_t}),   
u^{(k)}(x_1,t_{N_t}),\ldots,  u^{(k)}   (x_{N_x},t_{N_t})\big]^\mathsf{T}.
\end{aligned}
\end{align*}  
Then the discretization of the heat equation using SBP-SAT scheme on element $k$ leads to
\begin{align}\label{eq01}
\begin{aligned}
\big( \D{t} - \kappa^{(k)} {\D{x}}^2\big)\, \u^{(k)} 
	= \f^{(k)}-\S{w}{}-\S{e}{}- \S{0}{(k)} - \S{\Gamma}{(k,k-1)} - \S{\Gamma}{(k,k+1)},
\end{aligned}
\end{align}
where ${\D{x}}^2=\D{x}\D{x}$ and \(\f^{(k)}\) contains the source data.  The SAT terms \(\S{w}{},\, \S{e}{},\, \S{0}{(k)}\) and \(\S{\Gamma}{(k,k-1)}\) and \(\S{\Gamma}{(k,k+1)}\)  are defined to enforce boundary, initial, and interface conditions in the following way. The coupling between subdomains using these SATs is illustrated in Figure~\ref{fig:domaincoupling}.

To enforce the boundary conditions on the left and right ends of the domain, we add  $\S{w}{}$ and $\S{e}{}$ to the the first and last elements, respectively:
\begin{align}\label{eq01b}
\begin{aligned}
	\S{w}{} &=	\sigma_w\,\Pnorm{}^{-1}\,\R{w}^\mathsf{T}\,\PoneD{t}\,\bigl(\R{w}\u^{(1)}  -  \R{w}\h^{(1)}\bigr),\\[1ex]
	\S{e}{} &= \sigma_e\,\Pnorm{}^{-1}\,\R{e}^\mathsf{T}\,\PoneD{t}\,\bigl(\R{e}\u^{(K)} - \R{e}\g^{(K)}\bigr),
\end{aligned}
\end{align}
where \(\h^{(1)}\) and \(\g^{(K)}\) contain the boundary functions $h$ and $g$ evaluated at grid points at the first and last elements, and  \(\sigma_w\) and \(\sigma_e\) are positive penalty coefficients chosen to impose the boundary conditions weakly. These terms are set to be zero on all other elements. The SAT term for imposing the initial condition is  defined as 
{\setlength{\abovedisplayskip}{7pt}
 \setlength{\belowdisplayskip}{7pt}
\begin{align}\label{eq01t}
\begin{aligned}
	\S{0}{(k)}=\sigma_0\,\Pnorm{}^{-1}\,\R{s}^\mathsf{T}\,\PoneD{x}\,\bigl(\R{s}\u^{(k)} - \R{s}\q^{(k)}\bigr),\qquad\qquad {\tiny1 \le k \le K,}
\end{aligned}
\end{align}}
where \(\q^{(k)}\) contains the initial function $q$ evaluated at grid points on element $k$ and $\sigma_0$ is a positive constant.
The SAT terms for the interfaces are defined as
{\setlength{\abovedisplayskip}{7pt}
 \setlength{\belowdisplayskip}{7pt}
\begin{align}\label{eq01g}
\begin{aligned}
	&\S{\Gamma}{(k,k-1)} = \sigma_1
	\Pnorm{}^{-1}\,\R{w}^\mathsf{T}\,\PoneD{t}\,
	\bigl(\R{w}\u^{(k)} - \R{e}\u^{(k-1)}\bigr)\\
	&\quad + \sigma_2 \Pnorm{}^{-1}\,\R{w}^\mathsf{T}\,\PoneD{t}\,\bigl(\kappa^{(k)}\,\R{w}\,\D{x}\u^{(k)}
	- \kappa^{(k-1)}\, \R{e}\,\D{x}\u^{(k-1)}\bigr),\quad {\tiny 2 \le k \le K,}\\
    &\quad +\tau_1\,\kappa^{(k)}
	\Pnorm{}^{-1}\, \D{x}\Tr\R{w}^\mathsf{T}\,\PoneD{t}\,
	\bigl(\R{w}\u^{(k)} - \R{e}\u^{(k-1)}\bigr),
	\\&\\
	&\S{\Gamma}{(k,k+1)} =\sigma_3
	\Pnorm{}^{-1}\,\R{e}^\mathsf{T}\,\PoneD{t}\,
	\bigl(\R{e}\u^{(k)}-\R{w}\u^{(k+1)}\bigr)\\
	& \quad+\sigma_4 \Pnorm{}^{-1}\,\R{e}^\mathsf{T}\,\PoneD{t}\,\bigl(\kappa^{(k)}\,\R{e}\,\D{x}\u^{(k)}-\kappa^{(k+1)}\,\R{w}\,\D{x}\u^{(k+1)}
	\bigr),
	\quad {\tiny1 \le k \le K-1,}\\
     &\quad +\tau_2\,\kappa^{(k)}\Pnorm{}^{-1}\, \D{x}\Tr\R{e}^\mathsf{T}\,\PoneD{t}\,
	\bigl(\R{e}\u^{(k)}-\R{w}\u^{(k+1)}\bigr),
\end{aligned}
\end{align}}
where \(\sigma_1,\, \sigma_2,\, \sigma_3,\, \sigma_4, \, \tau_1 \) and \(\tau_2\) are  SAT coefficients. 

Notice that in the proposed SBP-SAT discretization, we use one temporal element in the time direction. Temporal refinement is performed by increasing the number of temporal collocation nodes. This choice keeps the energy analysis and the presentation transparent, and allows us to isolate the effect of the high-order space--time discretization. The SBP--SAT construction can also be generalized to multiple temporal elements by partitioning the time interval into time slabs and coupling neighbouring slabs weakly through temporal SAT terms. A discretization with multiple temporal elements  is relevant for large-scale space--time solvers, where the decomposition of the time interval provides a natural block structure for the global system and parallelization.

To construct the global linear system, let
{\setlength{\abovedisplayskip}{7pt}
 \setlength{\belowdisplayskip}{7pt} 
 \[\hat A^{(k)}=\Pnorm{}\D{t} - \kappa^{(k)} \Pnorm{}{\D{x}}^2 +\sigma_0\,\R{s}^\mathsf{T}\,\PoneD{x}\R{s}, \qquad \hat \f^{(k)}=\Pnorm{}\f^{(k)}+ \sigma_0 \,\R{s}^\mathsf{T}\,\PoneD{x}\R{s}\q^{(k)},\]}
where we have multiplied it by  $\Pnorm{}$ to avoid forming and applying  $\PoneD{x}{}^{-1}$ and $\PoneD{t}{}^{-1}$  in the SAT terms, improving efficiency and numerical robustness. Then define
\begin{align*}
	\begin{aligned}
		A^{(k)}&= \begin{cases}
			\hat A^{(k)}
			+\sigma_w \R{w}^\mathsf{T}\,\PoneD{t}\,\R{w}
			+\sigma_3\, \R{e}^\mathsf{T}\,\PoneD{t}\,\R{e}\\
			\quad\quad+\,\sigma_4\, \kappa^{(k)}\, \R{e}^\mathsf{T}\,\PoneD{t}\,\R{e}\,\D{x}
			+\tau_2\,\kappa^{(k)}{\D{x}}\Tr\R{e}^\mathsf{T}\,\PoneD{t}\,\R{e},
			&\; k=1,\\[1ex]
			\hat A^{(k)}
			+\sigma_1 \R{w}^\mathsf{T}\,\PoneD{t}\,\R{w}
			+\sigma_2 \,\kappa^{(k)}\, \R{w}^\mathsf{T}\,\PoneD{t}\,\R{w}\,\D{x}
			+\tau_1 \,\kappa^{(k)}{\D{x}}\Tr\R{w}^\mathsf{T}\,\PoneD{t}\,\R{w}
			\\
			\quad\quad
			+\, \sigma_3\, \R{e}^\mathsf{T}\,\PoneD{t}\,\R{e}
			+\sigma_4\, \kappa^{(k)}\, \R{e}^\mathsf{T}\,\PoneD{t}\,\R{e}\,\D{x}
			+\tau_2\,\kappa^{(k)}{\D{x}}\Tr\R{e}^\mathsf{T}\,\PoneD{t}\,\R{e},
			&\;2\le k\le K-1,\\[1ex]
			\hat A^{(k)}
			+\sigma_e \,\R{e}^\mathsf{T}\,\PoneD{t}\,\R{e}
			+\sigma_1 \R{w}^\mathsf{T}\,\PoneD{t}\,\R{w}\\
			\quad\quad+\,\sigma_2 \,\kappa^{(k)}\, \R{w}^\mathsf{T}\,\PoneD{t}\,\R{w}\,\D{x}
			+\tau_1 \,\kappa^{(k)}{\D{x}}\Tr\R{w}^\mathsf{T}\,\PoneD{t}\,\R{w},
			&\; k=K.
		\end{cases}
\end{aligned}
\end{align*}
and define
\begin{align*}
	\begin{aligned}      
		\b^{(k)}&=\begin{cases} 
			\hat \f^{(k)}+\sigma_w\,\R{w}^\mathsf{T}\,\PoneD{t}\,\R{w}\,\h^{(k)},&\; k=1,\\
			\hat \f^{(k)},&\; 2\le k\le K-1,\\
			\hat \f^{(k)}+\sigma_e\,\R{e}^\mathsf{T}\,\PoneD{t}\,\R{e}\,\g^{(k)},&\; k=K.\\
		\end{cases}
    \end{aligned}
\end{align*}
Denote
\begin{align*}
	\begin{aligned}        
		B^{(k)}&=-\sigma_3 \R{e}^\mathsf{T}\,\PoneD{t}\,\R{w}- \sigma_4\,  \kappa^{(k+1)}\R{e}^\mathsf{T}\,\PoneD{t} \R{w}\,\D{x}\,-\tau_2 \,\kappa^{(k)}{\D{x}}\Tr\R{e}^\mathsf{T}\,\PoneD{t}\,\R{w},\quad 1\le k\le K-1,\\
		C^{(k)}&=-\sigma_1 \R{w}^\mathsf{T}\,\PoneD{t}\,\R{e}- \sigma_2\,\, \kappa^{(k-1)} \R{w}^\mathsf{T}\,\PoneD{t} \R{e}\,\D{x}\,-\tau_1 \,\kappa^{(k)}{\D{x}}\Tr\R{w}^\mathsf{T}\,\PoneD{t}\,\R{e},\quad 2\le k\le K.	
	\end{aligned}
\end{align*}
The global linear system is given by	
\[
\begin{bmatrix}
	A^{(1)} & B^{(1)} & 0        & \cdots   & 0        \\[0.7ex]
	C^{(2)} & A^{(2)} & B^{(2)}  &          & \vdots   \\[0.7ex]
	0       & \ddots  & \ddots   & \ddots   & 0        \\[0.7ex]
	\vdots  &         & C^{(K-1)}& A^{(K-1)}& B^{(K-1)}\\[0.7ex]
	0       & \cdots  & 0        & C^{(K)}  & A^{(K)}
\end{bmatrix}
\;
\begin{bmatrix}
	\u^{(1)} \\[1.6ex]
	\u^{(2)} \\[1.6ex]
	\vdots  \\[1.6ex]
	\u^{(K)}
\end{bmatrix}
=\begin{bmatrix}
	\b^{(1)} \\[1.6ex]
	\b^{(2)} \\[1.6ex]
	\vdots  \\[1.6ex]
	\b^{(K)}
\end{bmatrix},
\]
or in residual form by
\begin{equation}\label{res}
	\calR(\u,  \vrho ) = \mathcal{A}(\vrho)\u - \b.
\end{equation}	
where $\vrho=[\rho^{(1)};\dots;\rho^{(K)}]\in[0,1]^K$ denotes the design vector.

\subsection{Conditional stability of proposed ST-SE SBP-SAT scheme}\label{sec:stab}

The following lemma is needed in the proof of the stability for the space-time SBP-SAT scheme.

\begin{lem}\label{lem:sbp-trace}(Discrete trace inequality)
For any grid vector $\z$,
\[\|\,\R{w}\z\,\|_{\PoneD{t}}^{2}\ \le\ \frac{1}{\hat p_{0}}\ \|\z\|^2_{\Pnorm{}},\qquad
\|\,\R{e}\z\,\|_{\PoneD{t}}^{2}\ \le\ \frac{1}{\hat p_{N_x}}\ \|\z\|^2_{\Pnorm{}},
\]
where \(\hat p_0= \e_w\Tr\, \PoneD{x}\e_w\), and \(\hat p_{N_x}= \e_e\Tr\, \PoneD{x}\e_e\).	
\end{lem}
\begin{proof} Let $\PoneD{x}=\diag(\hat p_0,\cdots,\hat p_{N_x})$ and $\PoneD{t}=\diag(\hat q_0,\cdots,\hat q_{N_t})$, so \(\hat p_0= \e_w\Tr\, \PoneD{x}\e_w\) and  \(\hat p_{N_x}= \e_e\Tr\, \PoneD{x}\e_e\).
First, note that the largest constant $c$ that can satisfy $c\,({\e_w}\Tr z)^2\le \z\Tr\,\PoneD{x}\z$ for all $\z\in\Rs^{N_x}$ is $c=\hat p_0$, since if we assume $\z=\alpha \e_w$, then 
\(\z\Tr\PoneD{x}\z=\alpha^2 {\e_w}\Tr\,\PoneD{x}\e_w=\alpha^2\hat p_0,\)
and 
\(({\e_w}\Tr z)^2=({\e_w}\Tr \alpha\, \e_w)({\e_w}\Tr \alpha \,\e_w)=\alpha^2.\)
Then substituting into $c\,({\e_w}\Tr z)^2\le \z\Tr\,\PoneD{x}\z$ yields $c\,\alpha^2\le \hat p_0\alpha^2$, forcing $c\le\hat p_0$, so the largest possible value is $c=\hat p_0$.
Now,
\begin{align*}
    \begin{aligned}
        \big\|\R{w}\z\big\|^2_{\PoneD{t}}&\le \sum_{j=1}^{N_t}\hat q_{j}(\R{w}\z)_j^2
        \le \sum_{j=1}^{N_t}\hat q_j \sum_{i=1}^{N_x}\frac{\hat p_i}{\hat p_0}z_i^2
       \le\frac{1}{\hat p_0}\sum_{j=1}^{N_t\,\times Nx}\Pnorm{jj} z_j^2=\frac{1}{\hat p_0} \|\z\|^2_{\Pnorm{}},
    \end{aligned}
\end{align*}
since $\Pnorm{}=\PoneD{t}\otimes\PoneD{x}$. The other inequality can be established in the same way. 
\end{proof}

\begin{thm}\label{thm:stab}
Assume the SBP-SAT scheme given by \eqref{eq01}--\eqref{eq01g} with zero source term. Let \(\hat p_0= \e_w\Tr\, \PoneD{x}\e_w\) and  \(\hat p_{N_x}= \e_{e}\Tr\, \PoneD{x}\e_{e}\). Then under the conditions
{\setlength{\abovedisplayskip}{7pt}
\setlength{\belowdisplayskip}{7pt}
\begin{align}\label{SATcond}
\begin{aligned}
&\sigma_0>\frac12,\qquad \sigma_w\ge\frac{\kappa_{\max}}{2\hat p_0},\qquad \sigma_e\ge\frac{\kappa_{\max}}{2\hat p_{N_x}},\qquad \sigma_1=\sigma_3\ge 0,\\
&\sigma_2=s,\quad \sigma_4=1+ s,\quad\tau_1=-(1+s),\quad \tau_2=-s,\qquad \forall s\in \Rs^+,
\end{aligned}
\end{align}}
the scheme is stable and following energy estimate holds:
{\setlength{\abovedisplayskip}{7pt}
 \setlength{\belowdisplayskip}{7pt}
\begin{align}\label{stab}
\begin{aligned}
&\sum_{k=1}^K\big\|\R{n}\u^{(k)}\big\|^2_{\PoneD{x}}
\le
\frac{\sigma_0^2}{2\sigma_0-1}
\sum_{k=1}^K\big\|\R{s}\q^{(k)}\big\|^2_{\PoneD{x}}
+\sigma_w\big\|\R{w}\h^{(1)}\big\|_{\PoneD{t}}^2
+\sigma_e\big\|\R{e}\g^{(K)}\big\|_{\PoneD{t}}^2 .
\end{aligned}
\end{align}}
\end{thm}

\begin{proof}
To derive a discrete energy estimate for the space–time scheme, we left–multiply \eqref{eq01} on each subdomain by $(\u^{(k)})\Tr \Pnorm{}$ and then sum over $k=1,\dots,K$. 
\begin{align}\label{ap01}
\begin{aligned}
\sum_{k=1}^K&{\u^{(k)}}\Tr \Pnorm{} \D{t} \u^{(k)}- \sum_{k=1}^K\kappa^{(k)} {\u^{(k)}}\Tr \Pnorm{}{\D{x}}^2 \u^{(k)}\\
    &\;=-\sum_{k=1}^K{\u^{(k)}}\Tr \Pnorm{}\S{0}{(k)}- 
    {\u^{(1)}}\Tr \Pnorm{}\S{w}{}-{\u^{(K)}}\Tr \Pnorm{}\S{e}{}
    -\sum_{k=2}^K{\u^{(k)}}\Tr \Pnorm{}\S{\Gamma}{(k,k-1)}-\sum_{k=1}^{K-1}{\u^{(k)}}\Tr \Pnorm{}\S{\Gamma}{(k,k+1)}.
\end{aligned}
\end{align}
By the SBP property,
\begin{align}\label{ap01a}
{\u^{(k)}}\Tr \Pnorm{}\D{t}\u^{(k)}=\tfrac12\,{\u^{(k)}}\Tr\E{t}\u^{(k)}
=\tfrac12\,(\R{n}\u^{(k)})\Tr\,\PoneD{x}\,(\R{n}\u^{(k)})-\tfrac12\,(\R{s}\u^{(k)})\Tr\,\PoneD{x}\,(\R{s}\u^{(k)}).
\end{align}
Moreover, since $ \Pnorm{}{\D{x}}^2=\E{x}\D{x}-{\D{x}}\Tr\Pnorm{}\D{x}$, we have
\begin{align*}
-{\u^{(k)}}\Tr \Pnorm{}{\D{x}}^2\u^{(k)}&=-{\u^{(k)}}\Tr\E{x}\D{x}\u^{(k)}
+{\u^{(k)}}\Tr {\D{x}}\Tr\Pnorm{}\D{x}\u^{(k)}\\
&=(\R{w}\u^{(k)})\Tr\,\PoneD{t}\,(\R{w}\D{x}\u^{(k)})-(\R{e}\u^{(k)})\Tr\,\PoneD{t}\,(\R{e}\D{x}\u^{(k)})+ (\D{x}\u^{(k)})\Tr\Pnorm{}(\D{x}\u^{(k)}).
\end{align*}
Substituting this and \eqref{ap01a} into \eqref{ap01} and taking sum over $k$ yields
\begin{align}\label{ap02}
\begin{aligned}
\frac12\sum_{k=1}^K\,&(\R{n}\u^{(k)})\Tr\,\PoneD{x}\,(\R{n}\u^{(k)})+\sum_{k=1}^K\kappa^{(k)}{\u^{(k)}}\Tr {\D{x}}\Tr\Pnorm{}\D{x}\u^{(k)}
\\=&
\sum_{k=1}^K\kappa^{(k)}(\R{e}\u^{(k)})\Tr\,\PoneD{t}\,(\R{e}\D{x}\u^{(k)})
-\sum_{k=1}^K\kappa^{(k)}(\R{w}\u^{(k)})\Tr\,\PoneD{t}\,(\R{w}\D{x}\u^{(k)})
\\&+ \frac12\sum_{k=1}^K\,(\R{s}\u^{(k)})\Tr\,\PoneD{x}\,(\R{s}\u^{(k)})
- {\u^{(1)}}\Tr \Pnorm{}\S{w}{}-{\u^{(K)}}\Tr \Pnorm{}\S{e}{}
\\&-\sum_{k=1}^K {\u^{(k)}}\Tr \Pnorm{}\S{0}{(k)}
- \sum_{k=2}^K{\u^{(k)}}\Tr \Pnorm{}\S{\Gamma}{(k,k-1)} -\sum_{k=1}^{K-1}{\u^{(k)}}\Tr \Pnorm{}\S{\Gamma}{(k,k+1)}.
\end{aligned}
\end{align}
It is straightforward to verify that
\begin{align}\label{ap03}
\begin{aligned}
	-{\u^{(1)}}\Tr\,\Pnorm{}\S{w}{}&= -\sigma_w\big\|\R{w}\u^{(1)}-\tfrac12\R{w}\h^{(1)}\big\|^2_{\PoneD{t}}
	\;+\;\tfrac{\sigma_w}{4}\,\big\|\R{w}\h^{(1)}\big\|_{\PoneD{t}}^2,\\
	-{\u^{(K)}}\Tr\,\Pnorm{}\S{e}{}&= -\sigma_e\big\|\R{e}\u^{(K)}-\tfrac12\R{e}\g^{(K)}\big\|^2_{\PoneD{t}}
	\;+\;\tfrac{\sigma_e}{4}\,\big\|\R{e}\g^{(K)}\big\|_{\PoneD{t}}^2,
\end{aligned}
\end{align}
and, for any $1\le k\le K$,
\begin{align}\label{ap04}
\begin{aligned}
\tfrac12\,(\R{s}\u^{(k)})\Tr&\,\PoneD{x}\,(\R{s}\u^{(k)})
-{\u^{(k)}}\Tr\,\Pnorm{}\S{0}{(k)}\\
&=
-\;\tfrac{2\sigma_0-1}{2}
\big\|\R{s}\u^{(k)}
-\tfrac{\sigma_0}{2\sigma_0-1}\R{s}\q^{(k)}\big\|^2_{\PoneD{x}}
+\;\tfrac{\sigma_0^2}{2(2\sigma_0-1)}
\big\|\R{s}\q^{(k)}\big\|^2_{\PoneD{x}}.
\end{aligned}
\end{align}
Substituting \eqref{ap03} and \eqref{ap04} into \eqref{ap02} gives
\begin{align}\label{ap05}
\begin{aligned}
\frac12\,&\sum_{k=1}^K\big\|\R{n}\u^{(k)}\big\|^2_{\PoneD{x}}
+\sum_{k=1}^K\kappa^{(k)}\big\|\D{x}\u^{(k)}\big\|^2_{\Pnorm{}}\\
=&
\sum_{k=1}^K\kappa^{(k)}(\R{e}\u^{(k)})\Tr\,\PoneD{t}\,(\R{e}\D{x}\u^{(k)})
-\sum_{k=1}^K\kappa^{(k)}(\R{w}\u^{(k)})\Tr\,\PoneD{t}\,(\R{w}\D{x}\u^{(k)})
\\&  -\sigma_w\big\|\R{w}\u^{(1)}-\tfrac12\R{w}\h^{(1)}\big\|^2_{\PoneD{t}}
+\;\frac{\sigma_w}{4}\,\big\|\R{w}\h^{(1)}\big\|_{\PoneD{t}}^2
 \\& -\sigma_e\big\|\R{e}\u^{(K)}-\tfrac12\R{e}\g^{(K)}\big\|^2_{\PoneD{t}}
 +\;\frac{\sigma_e}{4}\,\big\|\R{e}\g^{(K)}\big\|_{\PoneD{t}}^2
 \\&  -\;\frac{2\sigma_0-1}{2}\sum_{k=1}^K
 \big\|\R{s}\u^{(k)}
 -\tfrac{\sigma_0}{2\sigma_0-1}\R{s}\q^{(k)}\big\|^2_{\PoneD{x}}
 +\;\frac{\sigma_0^2}{2(2\sigma_0-1)}
 \sum_{k=1}^K\big\|\R{s}\q^{(k)}\big\|^2_{\PoneD{x}}
 \\&- \sum_{k=2}^K{\u^{(k)}}\Tr \Pnorm{}\S{\Gamma}{(k,k-1)}
 -\sum_{k=1}^{K-1}{\u^{(k)}}\Tr \Pnorm{}\S{\Gamma}{(k,k+1)}.
\end{aligned}
\end{align}
Then using Young's inequality, for any $\alpha>0$,
\begin{align}\label{ap06}
\begin{aligned}
-\kappa^{(1)}(\R{w}\u^{(1)})\Tr\,\PoneD{t}\,(\R{w}\D{x}\u^{(1)})\le \tfrac{\alpha \kappa^{(1)}}{2}\big\|\R{w}\u^{(1)}\big\|^2_{\PoneD{t}}+\tfrac{\kappa^{(1)}}{2\alpha}\big\|\R{w} \D{x}\u^{(1)}\big\|^2_{\PoneD{t}}
\end{aligned}
\end{align}
Notice that by using Lemma \eqref{lem:sbp-trace}
\begin{align}\label{ap07}
    \begin{aligned}
     \big\|\R{w} \D{x}\u^{(1)}\big\|^2_{\PoneD{t}}\le\tfrac{1}{\hat p_0} \big\|\D{x}\u^{(1)}\big\|^2_{\Pnorm{}},
    \end{aligned}
\end{align}
also
\begin{align*}
\begin{aligned}
&\big\|\R{w}\u^{(1)}\big\|^2_{\PoneD{t}}\le 2\big\|\R{w}\u^{(1)}-\tfrac12\R{w}\h^{(1)}\big\|^2_{\PoneD{t}}+\tfrac12\big\|\R{w}\h^{(1)}\big\|_{\PoneD{t}}^2.
\end{aligned}
\end{align*}
Then we choose $\alpha= \sigma_w/\kappa^{(1)}$ and substitute this and \eqref{ap07} in \eqref{ap06} to obtain
\begin{align}\label{ap09}
\begin{aligned}
-\kappa^{(1)}(\R{w}\u^{(1)})\Tr\,\PoneD{t}\,(\R{w}\D{x}\u^{(1)})&\le \sigma_w\big\|\R{w}\u^{(1)}-\tfrac12\R{w}\h^{(1)}\big\|^2_{\PoneD{t}}
\\&\quad+\tfrac{\sigma_w}{4}\big\|\R{w}\h^{(1)}\big\|_{\PoneD{t}}^2
+\tfrac{{\kappa^{(1)}}^2}{2\sigma_w\hat p_0} \big\|\D{x}\u^{(1)}\big\|^2_{\Pnorm{}}.
\end{aligned}
\end{align}
Doing the same for \(\kappa^{(K)}(\R{e}\u^{(K)})\Tr\,\PoneD{t}\,(\R{e}\D{x}\u^{(K)})\) 
and substituting the result together with \eqref{ap09} into \eqref{ap05}, we obtain
\begin{align}\label{ap10}
\begin{aligned}
\frac12\,\sum_{k=1}^K&\big\|\R{n}\u^{(k)}\big\|^2_{\PoneD{x}}+\sum_{k=2}^{K-1}\kappa^{(k)}\big\|\D{x}\u^{(k)}\big\|^2_{\Pnorm{}}
\\+&\,\kappa^{(1)}(1-\frac{\kappa^{(1)}}{2\sigma_w\hat p_0} )\big\|\D{x}\u^{(1)}\big\|^2_{\Pnorm{}}+\kappa^{(K)}(1-\frac{\kappa^{(K)}}{2\sigma_e\hat p_{N_x}} )\big\|\D{x}\u^{(K)}\big\|^2_{\Pnorm{}}\\
\le&
\sum_{k=1}^{K-1}\kappa^{(k)}(\R{e}\u^{(k)})\Tr\,\PoneD{t}\,(\R{e}\D{x}\u^{(k)})
-\sum_{k=1}^{K-1}\kappa^{(k+1)}(\R{w}\u^{(k+1)})\Tr\,\PoneD{t}\,(\R{w}\D{x}\u^{(k+1)})
\\&+\;\frac{\sigma_w}{2}\,\big\|\R{w}\h^{(1)}\big\|_{\PoneD{t}}^2
   +\;\frac{\sigma_e}{2}\,\big\|\R{e}\g^{(K)}\big\|_{\PoneD{t}}^2
   +\;\frac{\sigma_0^2}{2(2\sigma_0-1)}\sum_{k=1}^K\big\|\R{s}\q^{(k)}\big\|^2_{\PoneD{x}}
\\&-\;\sum_{k=2}^K{\u^{(k)}}\Tr \Pnorm{}\S{\Gamma}{(k,k-1)} - \sum_{k=1}^{K-1} {\u^{(k)}}\Tr \Pnorm{}\S{\Gamma}{(k,k+1)}.
\end{aligned}
\end{align}
So far we need to assume 
\begin{align}\label{ap11}
\begin{aligned}
\sigma_0>\frac12,\qquad\sigma_w\ge\frac{\kappa_{\max}}{2\hat p_0},\qquad \sigma_e\ge\frac{\kappa_{\max}}{2\hat p_{N_x}},
\end{aligned}
\end{align}
since $\kappa^{(1)},\, \kappa^{(K)}\le \kappa_{\max}$ by definition. It remains to show a stability estimate for SATs at the interface. To simplify calculation we introduce a vector 
\begin{align*}
\begin{aligned}
        \mat d_k=
\begin{bmatrix}
    \R{e}\u^{(k)}\\
    \R{w}\u^{(k+1)}\\
    \R{e}\D{x}\u^{(k)}\\
    \R{w}\D{x}\u^{(k+1)}
\end{bmatrix}.
\end{aligned}
\end{align*}
Collecting the terms at the interface
\(\Gamma_k=\partial\Omega_k\cap\partial\Omega_{k+1}\), for \(k=1,\ldots,K-1\),
we obtain
\begin{align*}
\begin{aligned}
&\kappa^{(k)}(\R{e}\u^{(k)})\Tr\PoneD{t}(\R{e}\D{x}\u^{(k)})
-\kappa^{(k+1)}(\R{w}\u^{(k+1)})\Tr\PoneD{t}(\R{w}\D{x}\u^{(k+1)})
\\
&\quad
-{\u^{(k)}}\Tr \Pnorm{}\S{\Gamma}{(k,k+1)}
-{\u^{(k+1)}}\Tr \Pnorm{}\S{\Gamma}{(k+1,k)}
=
{\mat d_k}\Tr(\mat M_k\otimes\PoneD{t})\mat d_k,
\end{aligned}
\end{align*}
where 
\begin{align*}
\begin{aligned}
\begingroup
\setlength{\arraycolsep}{4pt}
\renewcommand{\arraystretch}{1.25}
\mat{M}_k=\tfrac12\begin{bmatrix}
-2\sigma_3  & \sigma_1+\sigma_3 & \kappa^{(k)}(-\sigma_4-\tau_2+1) & \kappa^{(k+1)}(\sigma_4+\tau_1)\\
\sigma_1+\sigma_3 & -2\sigma_1 & \kappa^{(k)}(\sigma_2+\tau_2) & \kappa^{(k+1)}(-\sigma_2-\tau_1-1)\\
\kappa^{(k)}(-\sigma_4-\tau_2+1) & \kappa^{(k)}(\sigma_2+\tau_2) & 0 & 0\\
\kappa^{(k+1)}(\sigma_4+\tau_1) & \kappa^{(k+1)}(-\sigma_2-\tau_1-1) & 0 & 0
\end{bmatrix}.
\endgroup
\end{aligned}
\end{align*}
In order to have a stable scheme, we need to find conditions on SAT parameters such that $\mat M_k$ is negative semi-definite. For that we need to cancel out every column with zero diagonal entries.
This leads to $\sigma_4=1+\sigma_2$, resulting in a system of equations with a one-parameter family of solutions 
\begin{eqnarray}\label{ap13}
	\sigma_2=s,\quad \sigma_4=1+ s,\quad\tau_1=-(1+s),\quad \tau_2=-s,\qquad \forall s\in \Rs^+.
\end{eqnarray}
Then $\mat M_k$ reduces to
\begin{align*}
\begin{aligned}
\begingroup
\setlength{\arraycolsep}{4pt}
\renewcommand{\arraystretch}{1.25}
\mat{M}_k=\begin{bmatrix}
-\sigma_3  & \frac{\sigma_1+\sigma_3}2 &0& 0\\
\frac{\sigma_1+\sigma_3}2 & -\sigma_1 & 0& 0\\
0 & 0 & 0 & 0\\
0 & 0 & 0 & 0
\end{bmatrix},
\endgroup
\end{aligned}
\end{align*}
This matrix has only non-positive eigenvalues if \begin{equation}\label{ap14}
\sigma_1=\sigma_3\ge0.
\end{equation}
Since
\(\PoneD{t}\) is positive definite, \(\mat M_k\otimes\PoneD{t}\) also has
only non-positive eigenvalues. Hence
\[
\mat M_k\otimes\PoneD{t}\le0.
\]
Then if \eqref{ap11}, \eqref{ap13} and \eqref{ap14} hold, we obtain the following estimate from \eqref{ap10} for the given scheme:
\begin{align*}
\begin{aligned}
\frac12\,\sum_{k=1}^K&\big\|\R{n}\u^{(k)}\big\|^2_{\PoneD{x}}+\sum_{k=2}^{K-1}\kappa^{(k)}\big\|\D{x}\u^{(k)}\big\|^2_{\Pnorm{}}
\\&\quad+\kappa^{(1)}(1-\frac{\kappa^{(1)}}{2\sigma_w\hat p_0} )\big\|\D{x}\u^{(1)}\big\|^2_{\Pnorm{}}+\kappa^{(K)}(1-\frac{\kappa^{(K)}}{2\sigma_e\hat p_{N_x}} )\big\|\D{x}\u^{(K)}\big\|^2_{\Pnorm{}}\\
&\le\;\frac{\sigma_0^2}{2(2\sigma_0-1)} \sum_{k=1}^K\big\|\R{s}\q^{(k)}\big\|^2_{\PoneD{x}} 
+\;\frac{\sigma_w}{2}\,\big\|\R{w}\h^{(1)}\big\|_{\PoneD{t}}^2
+\;\frac{\sigma_e}{2}\,\big\|\R{e}\g^{(K)}\big\|_{\PoneD{t}}^2,
\end{aligned}
\end{align*}
Dropping the nonnegative terms on the left-hand side and multiplying by \(2\)
gives \eqref{stab}.
\end{proof}
In the proof, for simplicity we assume  a zero source term $\f^{(k)}\equiv 0$; the nonzero case follows by applying Young’s inequality to $(\u^{(k)})\Tr \Pnorm{} \f^{(k)}$ and absorbing the resulting $\|\u^{(k)}\|_{\Pnorm{}}^{2}$ term into the diffusion energy via a discrete Poincar\'e inequality.

Notice that for primal stability, the SAT coefficients admit lower bounds determined by the mesh resolution, specifically the endpoint quadrature weights $\hat p_0$ and $\hat p_{N_x}$ of the one-dimensional SBP norm $\PoneD{x}$. These bounds also scale with the material properties through the maximal thermal diffusivity $\kappa_{\max}$.

\subsection{Discrete topology optimization problem}\label{sec:TopOpt}
We assume the domain \(\Omega\) is partitioned into \(K\) design cells \(\{\Omega_{k}\}_{k=1}^K\)  as shown in Figure \ref{fig:domainElements} and $\vrho \in[0,1]^K$ denotes the design vector. Let $\u\in \Rs^N$ denote the vector of stacked temperature unknowns at all $N$ space-time nodes.   The state $\u$ and the design $ \vrho $ are related by the discrete residual
\(\calR(\u,  \vrho ) = \zero,
\) assembled from the space-time SBP-SAT scheme described in Section~\ref{sec:res}.

The discrete objective is the space-time integrated squared temperature
\( \Jac(\u,  \vrho ) = \u\Tr\,\Pnorm{}\,\u, \)
where  the weight matrix $\Pnorm{}$ is a symmetric positive definite matrix (SBP-quadrature matrix). Using the SBP-quadrature norm to approximate the continuous space--time functional keeps the discrete output functional compatible with the SBP--SAT primal and adjoint discretizations. Together with stability and dual consistency of the primal--adjoint pair, this compatibility provides the setting in which improved accuracy of functional estimates is expected~\cite{Hicken2011,Berg2013,Worku2022}, as explained in the following Section \ref{sec:adj}.

Instead of treating the discrete residual explicitly as a constraint, we use a reduced optimization formulation where the discrete residual is solved at each design iteration \cite{Sigmund2004}. The discrete form  of the topology optimization problem is then given by
{\setlength{\abovedisplayskip}{7pt}
 \setlength{\belowdisplayskip}{7pt}
\begin{align}\label{eq110}
	\begin{aligned}
            & \min_{\vrho\in[0,1]^K } . \quad   \Jac(\u,  \vrho ) = \u\Tr\,\Pnorm{}\,\u \qquad \mbox{where}\; \u\; \mbox{satisfies}\;\calR(\u,  \vrho ) = \zero,\\
			& \;\;\text{s.t.} \quad  \sum_{k=1}^{K}\rho^{(k)}V_x^{(k)} \le V^*, 
	\end{aligned}
\end{align}}
where
\(V_x^{(k)}=b_k-a_k.\)
 
At each iteration, we first solve the forward problem $\calR(\u,\vrho)=\zero$ for $\u$ with the current design $\vrho$; followed by solving the adjoint equation  for adjoint vector variable $\v = [\v^{(1)};\cdots;\v^{(K)}]$
\[\displaystyle \Bigl(\dfrac{\partial\calR}{\partial \u}\Bigr)^{{\mathsf{T}}}\!\v
			=\dfrac{\partial\Jac}{\partial \u}.\]			      
Then we compute the design sensitivity 
\[\displaystyle \nabla_{\rho}\,\Jac
			= -\,\v^{{\mathsf{T}}}\,\frac{\partial\calR}{\partial \vrho}.\]
Finally, we update the design $\vrho$ using the method of moving asymptotes (MMA) \cite{Svanberg1987},  subject to the volume constraint. MMA is a widely adopted method for density-based topology optimization. 
Thus for brevity, we omit algorithmic details and refer the reader to \cite{Svanberg1987} for the method and its application to topology optimization.
Each iteration requires one forward and one adjoint solve, followed by computation of the sensitivity in all space–time elements. These steps are repeated until the desired convergence is achieved.

\subsection{Dual-consistent SBP-SAT scheme}\label{sec:adj}
In the following, we use the following notation. Let \(u\) denote a given smooth function, and let \(\u_{\rm ex}\) be the vector obtained by evaluating \(u\) at the space--time collocation nodes. For simplicity, we assume that the polynomial orders in time and space are refined together and write \(N:=N_t=N_x\). Let \(u_{\rm N}\) denote the ST--SE approximation in the space--time polynomial space, and let \(\u_{\rm N}\) denote its vector of nodal values at the space--time collocation nodes. Let \(\calL\) be the linear differential operator of the transient heat equation with piecewise-constant diffusion coefficient and  \(J(u)\) be a given output functional.  The adjoint differential operator is given by
\[
    \calL^{\,\!*}v = -v_t-(\kappa v_x)_x = g := \frac{\partial J}{\partial u},
\]
together with the appropriate adjoint boundary and terminal conditions~\cite{Lanczos1961}. For example, when the primal PDE imposes an initial condition at \(t=0\) and inhomogeneous Dirichlet boundary conditions, the adjoint PDE imposes a terminal condition at \(t=T\) and homogeneous Dirichlet boundary conditions; see Appendix~\ref{App:AdjSol}. Finally, let \(L_{\rm N}\) and \(L_{\rm N}^{\!*}\) be discrete operators approximating \(\calL\) and \(\calL^{*}\), respectively.

A discretization is \emph{dual-consistent} or \emph{adjoint-consistent} if it yields a discrete dual problem that is a consistent discretization of the continuous dual PDE~\cite{Hartmann2007}. 
In the spectral setting considered here, we assess this consistency by the exponential decay of the adjoint residual
\[
\|L_{\rm N}^{\!*}\v_{\rm ex}-\g_{\rm ex}\|_{\Pnorm{}}
=
\mathcal{O}\!\left(e^{-\beta_v N}\right),
\qquad \beta_v>0,
\quad \text{as } N\to\infty,
\]
where \(\v_{\rm ex}\) and \(\g_{\rm ex}\) denote the exact adjoint solution and adjoint source evaluated at the space--time collocation nodes.

In general, primal consistency of a discrete operator alone does not imply consistency of the discrete dual problem, and accurate / high-order convergence of target functionals typically requires a dual-consistent discretization. For algebraic convergence rates, \emph{superconvergence functional estimate} refers to the situation in which the functional error converges at a higher order than the solution error~\cite{Niles2000,Giles2002}. 
In the present spectral setting, we examine whether the functional error retains spectral decay under the dual-consistent discretization. More precisely, the solution approximation is said to converge spectrally if
\[
\|\u_{\rm N}-\u_{\rm ex}\|_{\Pnorm{}}
=
\mathcal{O}\!\left(e^{-\beta_u N}\right),
\qquad \beta_u>0,
\quad \text{as } N\to\infty.
\]
Similarly, the functional estimate is said to converge spectrally if
\[
|J(u_{\rm N})-J(u)|
=
\mathcal{O}\!\left(e^{-\beta_J N}\right),
\qquad \beta_J>0,
\quad \text{as } N\to\infty.
\]
If the observed decay rate of the functional error is larger than that of the solution error, i.e. \(\beta_J>\beta_u\), then the functional estimate exhibits faster spectral convergence.

In the proposed stable space--time formulation, the temporal SATs are chosen so that the transpose of the primal residual Jacobian consistently imposes the corresponding terminal condition for the adjoint problem. The primal initial condition is imposed weakly by
\[
\begin{aligned}
	\S{0}{(k)}=\sigma_0\,\Pnorm{}^{-1}\,\R{s}^\mathsf{T}\,\PoneD{x}\,\bigl(\R{s}\u^{(k)} - \R{s}\q^{(k)}\bigr),\qquad\qquad {\tiny1 \le k \le K,}
\end{aligned}
\]
while the adjoint terminal condition is imposed through
\[
\begin{aligned}
	\S{f}{(k)}=\sigma_f\,\Pnorm{}^{-1}\,\R{n}^\mathsf{T}\,\PoneD{x}\,\bigl(\R{n}\v^{(k)} - \zero\bigr),\qquad\qquad {\tiny1 \le k \le K.}
\end{aligned}
\]
In the reduced optimization formulation, the discrete adjoint is written as
\[
\Bigl(\dfrac{\partial \calR}{\partial \u}\Bigr)^{\mathsf{T}}\!\v
= \dfrac{\partial \Jac}{\partial \u}
\quad\Rightarrow\quad
\calA^{\mathsf{T}}\v = 2\,\Pnorm{}\,\u,
\]
where $\v = [\v^{(1)};\cdots;\v^{(K)}]$ is the  vector obtained by stacking the subdomain adjoint variables. Equivalently, let \(\hat\calA\) denote the unweighted SBP--SAT operator, \(\calA:=\Pnorm{}\hat\calA\), then the adjoint equation  \(\calA^{\mathsf{T}}\v = 2\,\Pnorm{}\,\u\) can be written as
\[
\hat\calA^\dagger\v=2\u,
\]
where \(\hat\calA^\dagger=\Pnorm{}^{-1}\hat\calA^\mathsf{T}\Pnorm{}\) is the adjoint of the associated matrix with respect to the SBP inner product.

Dual-consistent discretizations and superconvergence functional estimate  have been studied for Galerkin and finite-difference SBP--SAT methods. 
For Galerkin and DG methods, adjoint-error-correction techniques have been developed to improve the accuracy of output functionals~\cite{Cockburn2007}. 
For finite-difference SBP discretizations, conditions for functional superconvergence have been derived for both steady and time-dependent problems~\cite{Hicken2011,Berg2012}.
For diffusion problems, related results have been obtained for narrow-stencil second-derivative generalized SBP operators~\cite{Worku2022}. 
These works indicate that, under suitable regularity, stability, and dual-consistency assumptions, discretized output functionals can exhibit improved accuracy compared with the underlying state approximation. 
In the present topology optimization setting, however, the diffusivity is generally piecewise constant  and discontinuous because it is induced by the element-wise design field \(\rho\). 
Thus, such functional-accuracy results can only be expected when material interfaces are adequately resolved, for example through mesh refinement, or when the design field is sufficiently regularized, for example by filtering~\cite{Sigmund2004}. 
Thus, while these works motivate the present dual-consistent construction, a complete functional-error analysis for the proposed tensor-product space--time spectral-element SBP--SAT discretization in multi-material domains is left for future work.

\subsection{Calculating design sensitivity}\label{sec:sen}
The design sensitivity is given by
\[\nabla_{\rho}\Jac= -\,\v^{\mathsf{T}}\,\frac{\partial \calR}{\partial \vrho},\]
where
\[
\nabla_{\rho}\Jac =\bigg[
\frac{\partial \Jac}{\partial \rho_1} ;\;
\frac{\partial \Jac}{\partial \rho_2} ;\;
\cdots ;\;
\frac{\partial \Jac}{\partial \rho_K}
\bigg].
\]
In particular, for each component $\rho_k$,
\[
\frac{\partial \Jac}{\partial \rho_k}
= -\,\sum_{j=1}^{K}
\bigg(\frac{\partial \calR_j}{\partial \rho_k}\bigg)^{\mathsf{T}} \v^{(j)},
\qquad 1 \le k \le K,
\]
where $\calR_j$ denotes the $j$-th block row of $\calR$. 
Notice that for the first element
\begin{align*}
\begin{aligned}
\frac{\partial \calR}{\partial \rho_1}\;&=\;
\begin{bmatrix}
	\frac{\partial A^{(1)}}{\partial \rho_1} & \frac{\partial B^{(1)}}{\partial \rho_1} & 0 & \cdots & 0 \\
	\frac{\partial C^{(2)}}{\partial \rho_1} & \frac{\partial A^{(2)}}{\partial \rho_1} & \frac{\partial B^{(2)}}{\partial \rho_1} &  \vdots      & \vdots \\
	0 & \ddots & \ddots & \ddots & 0 \\
	\vdots & \ddots & \frac{\partial C^{(K-1)}}{\partial \rho_1} & \frac{\partial A^{(K-1)}}{\partial \rho_1} & \frac{\partial B^{(K-1)}}{\partial \rho_1} \\
	0 & \cdots & 0 & \frac{\partial C^{(K)}}{\partial \rho_1} & \frac{\partial A^{(K)}}{\partial \rho_1}
\end{bmatrix}
\begin{bmatrix}
	\u^{(1)} \\[1.8ex]
	\u^{(2)} \\[1.8ex]
	\vdots \\[1.8ex]
	\u^{(K)}
\end{bmatrix}
-
\begin{bmatrix}
	\frac{\partial \,b^{(1)}}{\partial \rho_1} \\[1.2ex]
	\frac{\partial \,b^{(2)}}{\partial \rho_1} \\[1.2ex]
	\vdots \\[1.2ex]
	\frac{\partial \,b^{(K)}}{\partial \rho_1}
\end{bmatrix}\\
&=\;
\begin{bmatrix}
	\frac{\partial A^{(1)}}{\partial \rho_1} & \frac{\partial B^{(1)}}{\partial \rho_1}  & \cdots & 0 \\[1.2ex]
	\frac{\partial C^{(2)}}{\partial \rho_1} & 0  & \cdots & \vdots \\[1.2ex]
    0&0&\cdots&0\\
	\vdots & \vdots  & \ddots & \vdots \\
	0 & \cdots  &\cdots & 0
\end{bmatrix}
\begin{bmatrix}
	\u^{(1)} \\[1.2ex]
	\u^{(2)} \\[1.2ex]
	\vdots \\[1.2ex]
	\u^{(K)}
\end{bmatrix},
\end{aligned}
\end{align*}
then 
\begin{align*}
\begin{aligned}
    \dfrac{\partial \Jac}{\partial \rho_1}=
    \;-\left(\;
\begin{bmatrix}
	\frac{\partial A^{(1)}}{\partial \rho_1} & \frac{\partial B^{(1)}}{\partial \rho_1}    \\[1.2ex]
	\frac{\partial C^{(2)}}{\partial \rho_1} & 0  
\end{bmatrix}
\begin{bmatrix}
	\u^{(1)} \\[2.8ex]
	\u^{(2)}
\end{bmatrix}\right)^{{\mathsf{T}}}\cdot\,
\begin{bmatrix}
	\v^{(1)} \\[1.8ex]
	\v^{(2)} 
\end{bmatrix},
\end{aligned}
\end{align*}
since only the matrices for a given element depend on the design variable of that same element.
Similarly for a given internal element
\begin{align*}
\begin{aligned}
    \dfrac{\partial \Jac}{\partial \rho_k}=
    \;-\left(\;
\begin{bmatrix}
	0&\frac{\partial B^{(k-1)}}{\partial \rho_k} & 0   \\[1.2ex]
	\frac{\partial C^{(k)}}{\partial \rho_k} &\frac{\partial A^{(k)}}{\partial \rho_k} & \frac{\partial B^{(k)}}{\partial \rho_k}   \\[1.2ex]
	0&\frac{\partial C^{(k+1)}}{\partial \rho_k} & 0  
\end{bmatrix}
\begin{bmatrix}
	\u^{(k-1)} \\[2ex]
	\u^{(k)}   \\[2ex]
	\u^{(k+1)} 
\end{bmatrix}\right)^{{\mathsf{T}}}\cdot\,
\begin{bmatrix}
	\v^{(k-1)} \\[1.6ex]
	\v^{(k)}   \\[1.6ex]
	\v^{(k+1)} 
\end{bmatrix},\qquad 2\le k\le K-1,
\end{aligned}
\end{align*}
and for the final element
\begin{align*}
\begin{aligned}
    \dfrac{\partial \Jac}{\partial \rho_K}=
    \;-\left(\;
\begin{bmatrix}
	0&\frac{\partial B^{(K-1)}}{\partial \rho_K}    \\[1.2ex]
	\frac{\partial C^{(K)}}{\partial \rho_K} &\frac{\partial A^{(K)}}{\partial \rho_K}  
\end{bmatrix}
\begin{bmatrix}
\u^{(K-1)}\\[2.8ex]
\u^{(K)} 
\end{bmatrix}\right)^{{\mathsf{T}}}\cdot\,
\begin{bmatrix}
	\v^{(K-1)} \\[1.8ex]
	\v^{(K)} 
\end{bmatrix}.
\end{aligned}
\end{align*}

\section{Computational performance}\label{sec:results}
In the following we present numerical results to evaluate the performance and effectiveness of the proposed method. Throughout this paper, we compare the timing and convergence rates of two all-at-once forward solvers. The following abbreviations are used throughout the numerical results:
\begin{enumerate}
\item [(i)] 
BE--FE--AAO: backward Euler in time with linear finite element in space, solved as a single all-at-once linear system;

\item [(ii)]  ST--SE: the proposed space-time spectral element SBP-SAT formulation.
\end{enumerate}

The numerical results in Section~\ref{sec:num:exp} demonstrate the spectral convergence behavior of the proposed method in comparison with a low-order method. In Sections~\ref{sec:num:verify} and ~\ref{sec:num:dual} we provide numerical experiments that validate the optimal designs and dual consistency through comparisons with the independently computed reference solution. 
In Section~\ref{sec:num:time}, we report time-to-solution (wall-clock time) and provide cost-of-accuracy (memory requirement) curves by plotting the relative change in objective error versus runtime. All experiments were run in MATLAB~R2024b on a Windows~11 laptop with an AMD Ryzen~7~8840U with 8~cores at 3.30\,GHz and 32\,GB of RAM. 

\subsection{Exponential convergence of the 1D forward solver}\label{sec:num:exp}
\begin{figure}
  \centering
  \includegraphics[width=0.55\linewidth]{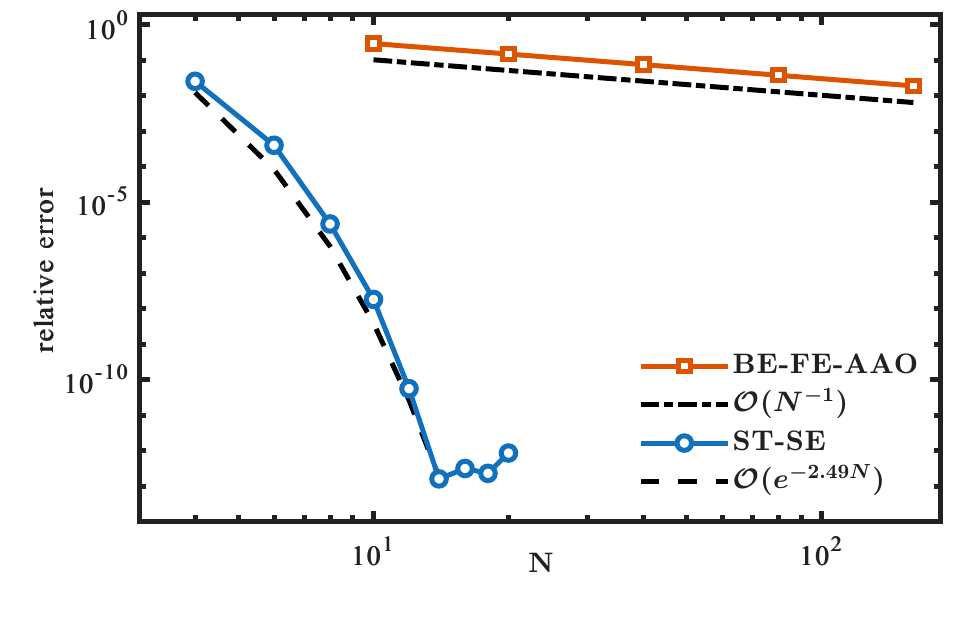}
  \caption{1D Forward-solver error comparison between BE--FE--AAO and ST--SE 
  methods, $N:=N_t=N_x$.}
  \label{fig:int}
\end{figure}

Figure~\ref{fig:int} compares the convergence of ST--SE and BE--FE--AAO for the one-dimensional heat equation on the space--time domain \([-2,1]\times[0,1]\). 
The test uses \(K=10\) space-time elements with different constant material properties, so that the material heterogeneity is resolved by element boundaries,
\[
\rho^{(k)} =
0.5 + 0.3\,\sin(2\pi\theta_k)
      + 0.15\,\cos(6\pi\theta_k),
\qquad
\theta_k=\frac{k-\frac12}{K},
\]
for $1\le k\le K$ and  we set \(p=3\) in~\eqref{eq0} with $\kappa_{\min}=0.05$ and $\kappa_{\max}=1$.
Both methods are compared using the same prescribed
heterogeneous conductivities and the same piecewise-defined manufactured
solution. In the \(k\)-th subdomain \([a_k,b_k]\), the solution is given by
\[
u^{(k)}(x,t)
=
(-1)^{k-1}
\frac{L_k}{\kappa^{(k)}}
\sin\!\left(\frac{\pi(x-a_k)}{L_k}\right)
e^{-\pi t},\quad L_k=b_k-a_k,
\]
for $1\le k\le K$ and $t\in[0,1]$. This construction satisfies the
homogeneous Dirichlet conditions at the external boundaries and ensures
continuity of both the solution and the conductive flux across the
subdomain interfaces. The SAT parameters in the ST--SE scheme are chosen to satisfy the stability conditions in Theorem~\ref{thm:stab}. In particular, in~\eqref{SATcond}, we set throughout all numerical examples
\[\sigma_0=1,\qquad \sigma_1=\sigma_3=1,\qquad \sigma_2=1,\qquad \sigma_4=2,\qquad \tau_1=-2,\qquad \tau_2=-1. \] The boundary SAT parameters $\sigma_w$ and $\sigma_e$ are chosen to satisfy the stability bounds determined by the endpoint weights of the one-dimensional SBP norm matrix $\PoneD{x}$ and the maximum conductivity $\kappa_{\max}$.

The spectral convergence of ST--SE is observed, with the error decreasing almost exponentially as the number of collocation nodes, $N:=N_t=N_x$, increases and then remaining around $10^{-14}$, where floating-point round-off errors dominate. In contrast, BE--FE--AAO exhibits the expected algebraic convergence associated with the low-order backward Euler and linear finite-element discretization. 

For the classical spectral element methods, spectral decay of the error is typically observed provided that the solution is a smooth function. This result numerically demonstrates that the space--time spectral element SBP--SAT formulation retains spectral convergence for smooth manufactured solutions in a heterogeneous domain with material interfaces resolved by element boundaries, whereas the low-order all-at-once BE--FE discretization converges at an algebraic rate.

\subsection{Verifying the optimal design in the numerical solver}\label{sec:num:verify}
We consider a simplified model topology optimization problem subject to the 1D transient heat equation given by
\begin{align}\label{primal}
\begin{aligned}
& u _t - (\kappa(x)  u _x)_x = f, &\qquad& (x,t)\in [0,1]\times[0,T],\\
& u (0,t) = 0,\quad  u (1,t)= u _R,&\qquad& t\in[0,T],\\
&u (x,0) =  u _0(x),&\qquad&x\in[0,1],
\end{aligned}
\end{align}
where  $ u _R$ and $f$ are given positive constants. Suppose there is an interface at \(\xi\in(0,1)\) and $\kappa(x)$ is a piecewise positive function given by,
\begin{align}\label{k2domain}
\begin{aligned}
\kappa(x)=
\begin{cases}
\kappa_1,& 0\le x\le\xi,\\
\kappa_2,& \xi<x\le1,
\end{cases}
\qquad \kappa_1,\,\kappa_2>0,
\end{aligned}
\end{align}
rendering the problem with just two design variables.
The interface conditions at $x=\xi$ are given by
\[
 u (\xi^-,t)= u (\xi^+,t), \qquad
\kappa_1\, u _x(\xi^-,t)=\kappa_2\, u _x(\xi^+,t).
\]
A manufactured solution consistent with this design is derived in Appendix~\ref{app:twoDomSol}, from which we obtain the initial condition \( u (x,0)= u _0(x)\).

We consider the objective functional given by
\[J( u )=\int_{0}^{T}\!\!\int_{0}^{1}  u (x,t)^2\,\mathrm{d}x\,\mathrm{d}t, \]
where $T=0.5$ and constrain the material by \(\kappa_1+\kappa_2=0.75\). We also   
assume that $\xi=0.5$ so that the domain is split into two equal-length subdomains. Then the optimization problem is reduced to
the one–dimensional search for \(\kappa_2=0.75-\kappa_1\), which is solved by the MATLAB-function \texttt{fminbnd}  with a tolerance $10^{-8}$, yielding an approximate analytical
minimizer 
\[\kappa_1=0.4341,\qquad \kappa_2= 0.3158,\] where the optimal value for the cost function is  \(\Jac=5.1730.\)  We refer to the corresponding solver as \texttt{MMSHeat}. The optimal allocation \(\kappa_1>\kappa_2\) is physically intuitive: lowering the space–time \(L^2\)-temperature requires evacuating heat toward the cold boundary at \(x=0\).

To make a comparison with this (approximate) analytical solution, we run the
topology optimization solver (which we denote by \texttt{MMAHeat}) with the same settings. MMA iterations were terminated when the maximum change in the design variables between successive iterates fell below $10^{-8}$, or after at most 100 iterations. To match the design variables with the $\kappa_1$ and $\kappa_2$ in analytical solution study, we remove the penalization effects in~\eqref{eq0}, setting \(p=1\)  (linear material interpolation), \(\kappa_{\max} = 1\) and \(\kappa_{\min} = 0\). We enforce the same material budget used in the analytical solution by choosing \(V^*=0.375\) in the
volume constraint  inequality ~\eqref{eq00} since it represents the mean–conductivity (due to being weighted by the element volume); on two equal elements this is equivalent to \(\kappa_1+\kappa_2=0.75\).

\begin{figure}
	\centering
    \includegraphics[width=0.48\linewidth]{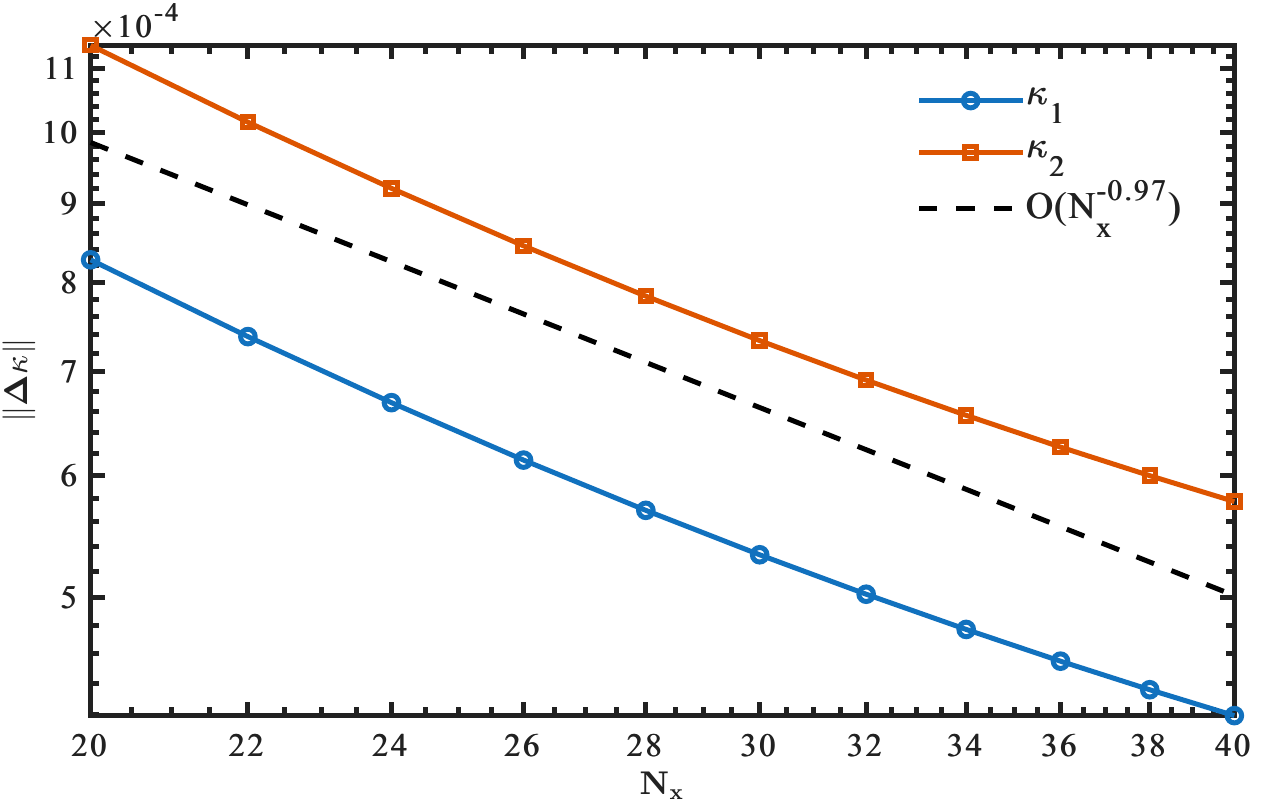}
    \includegraphics[width=0.48\linewidth]{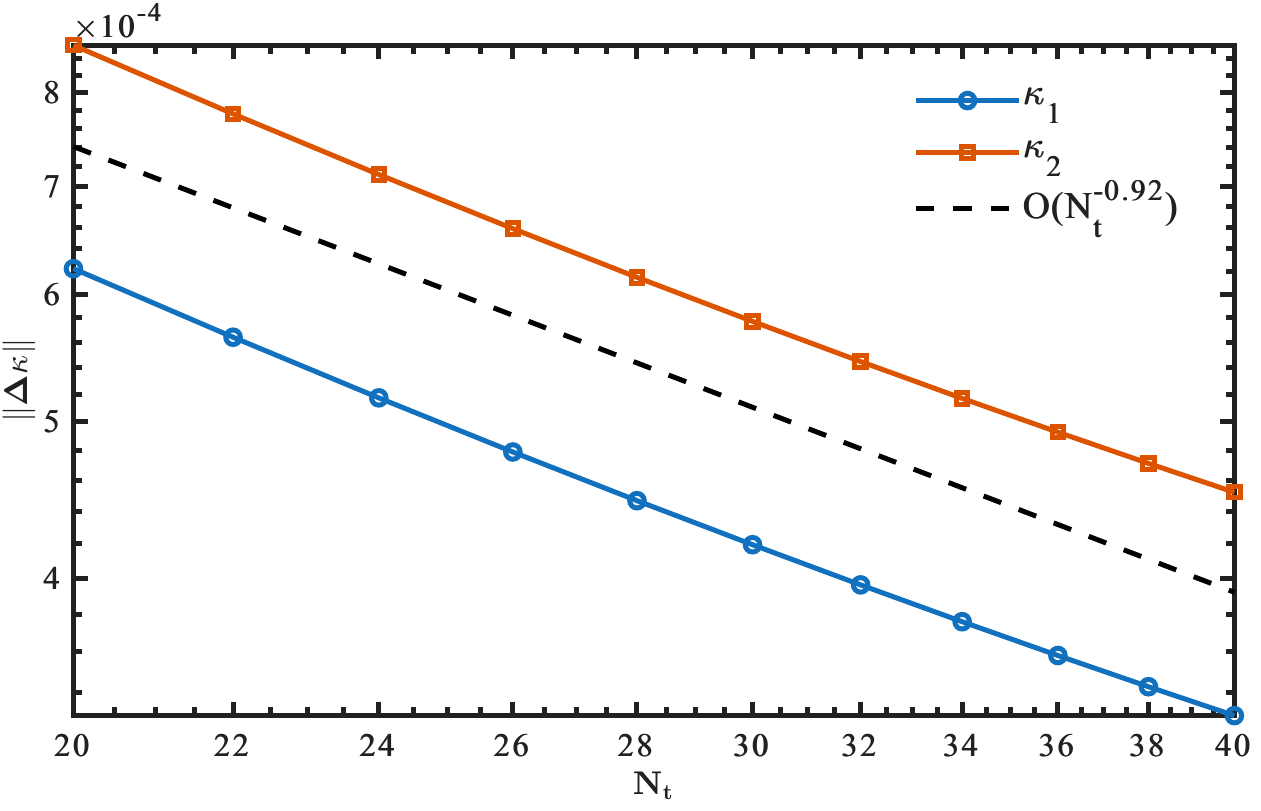}
    \caption{Relative error of the \texttt{MMAHeat} optimal design (using \texttt{MMSHeat} as the reference). Left: spatial refinement by increasing \(N_x\) with \(N_t=30\); right: temporal refinement by increasing \(N_t\) with \(N_x=40\). 
    }
    \label{fig:TwoDom1}
\end{figure}

Figure~\ref{fig:TwoDom1} reports the relative error of the \texttt{MMAHeat} optimal design, using \texttt{MMSHeat} as the reference (treated as exact in this test). The error decreases steadily under mesh refinement in one direction, while the other direction is held fixed. The observed linear convergence is due to comparing two numerical results from two optimization algorithms that they are converging linearly.

\subsection{Verifying dual consistency and functional accuracy}\label{sec:num:dual}

The purpose of the following experiment is to verify, for a fixed heterogeneous design, that the transpose of the proposed discrete primal SBP--SAT operator is a consistent discretization of the continuous adjoint equation given by   
\begin{align}\label{adjoint}
\begin{aligned}
&\calL^{\,\!*} v := -\,v_t - (\kappa v_x)_x = 2u, 
&\qquad &(x,t)\in [0,1]\times[0,T],\\
&v(0,t)=v(1,t)=0,
&\qquad &t\in[0,T],\\
&v(x,T)=0, 
&\qquad &x\in[0,1],
\end{aligned}
\end{align}
subject to the interface conditions at \(x=\xi\in[0,1]\),
\[
v(\xi^-,t)=v(\xi^+,t),\qquad 
\kappa_1 v_x(\xi^-,t)=\kappa_2 v_x(\xi^+,t),
\]
and where \(\kappa(x)\) is the piecewise positive function given by \eqref{k2domain}. 
For the adjoint consistency test, rather than using \(2u\) as the adjoint source, we prescribe the manufactured adjoint solution
\[
    v(x,t)=\sin^2(2\pi x)\sin^2(\pi t/T),
    \qquad t\in[0,T],
\]
and define the manufactured adjoint source by
\(g:=-v_t-(\kappa v_x)_x.\)
The errors then are defined by
\[
E_u =
\frac{\|\u_{\rm N}-\u_{\rm ex}\|_{\Pnorm{}}}
     {\|\u_{\rm ex}\|_{\Pnorm{}}},
\qquad
E_J =
\frac{
\left|\Jac_{\rm N}(\u_{\rm N},\vrho)-J_{\rm ex}\right|
}
{|J_{\rm ex}|},
\qquad
E_{\rm adj} =
\frac{
\left\|
\calA(\vrho)^{\mathsf T}\v_{\rm ex}
-
\Pnorm{} \g_{\rm ex}
\right\|_{\Pnorm{}^{-1}}
}{
\| \g_{\rm ex}\|_{\Pnorm{}}
}.
\]
Here \(J_{\rm ex}\) is the scalar reference value of the continuous functional, computed using high-order quadrature. 
The vectors \(\u_{\rm ex}\), \(\v_{\rm ex}\), and \(\g_{\rm ex}\) denote the nodal values of the manufactured primal solution of problem~\eqref{primal}, the manufactured adjoint solution, and the manufactured adjoint source, respectively, evaluated at the space--time LGL nodes.

\begin{figure}
	\centering
    \includegraphics[width=0.55\linewidth]{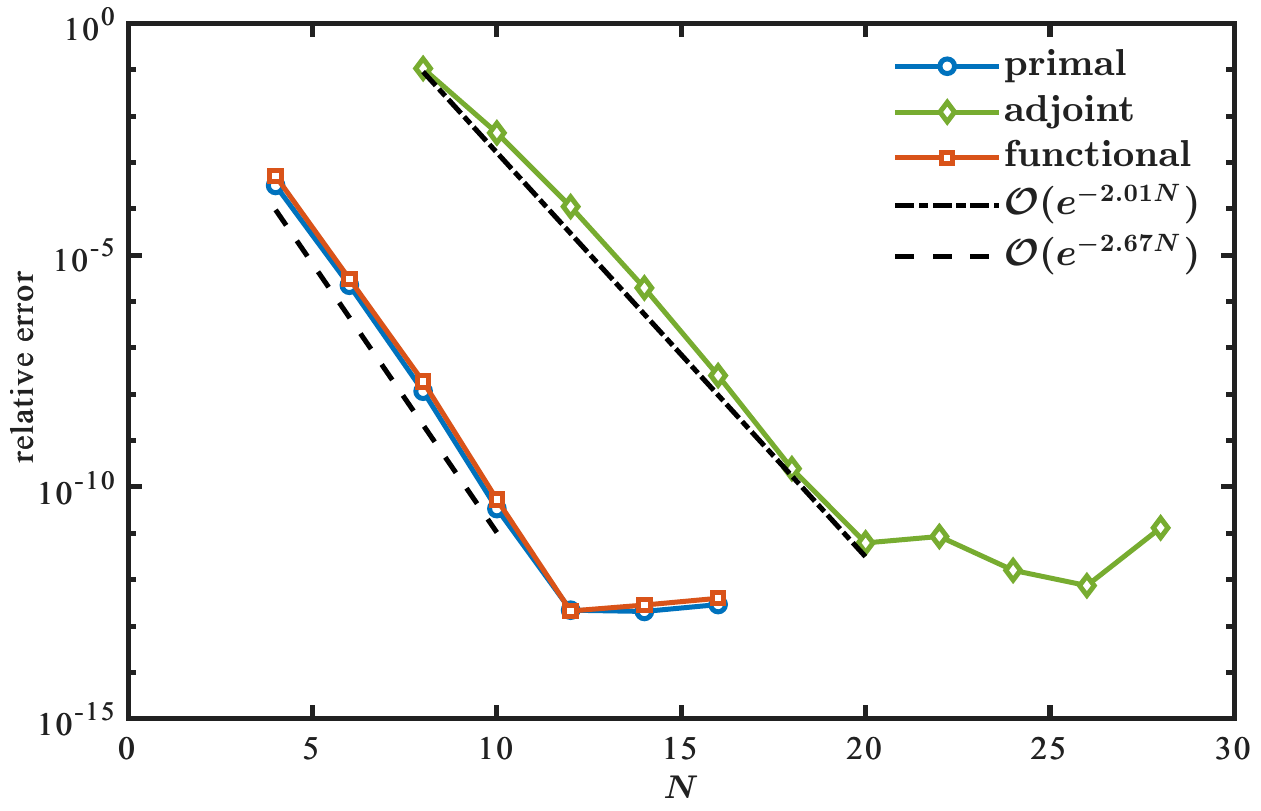}
	\caption{Fixed-design  primal--adjoint consistency test for ST--SE scheme. The primal state error, quadratic functional error, and residual-based adjoint consistency error decay exponentially, \(N := N_t = N_x.\)}
	\label{fig:dual}
\end{figure}

The left panel of Figure~\ref{fig:dual} reports a fixed-design convergence rate test for the primal, adjoint and functional discretizations under space-time refinement where $T=0.5$. In this test, the design vector \(\vrho\), and hence the material coefficient \(\kappa(\vrho)\), is prescribed and kept fixed; no optimization is performed. The polynomial order \(N := N_x = N_t\) is then increased. The numerical experiment shows that the primal error, functional error, and residual-based adjoint consistency error all decay exponentially with \(N\) until round-off effects dominate.  This verifies, for this smooth fixed-design test, that the primal space--time solver is spectrally accurate and that the transpose operator \(\calA(\vrho)^{\mathsf T}\) is consistent with the manufactured continuous adjoint equation.

\begin{figure}
	\centering
    \includegraphics[width=0.55\linewidth]{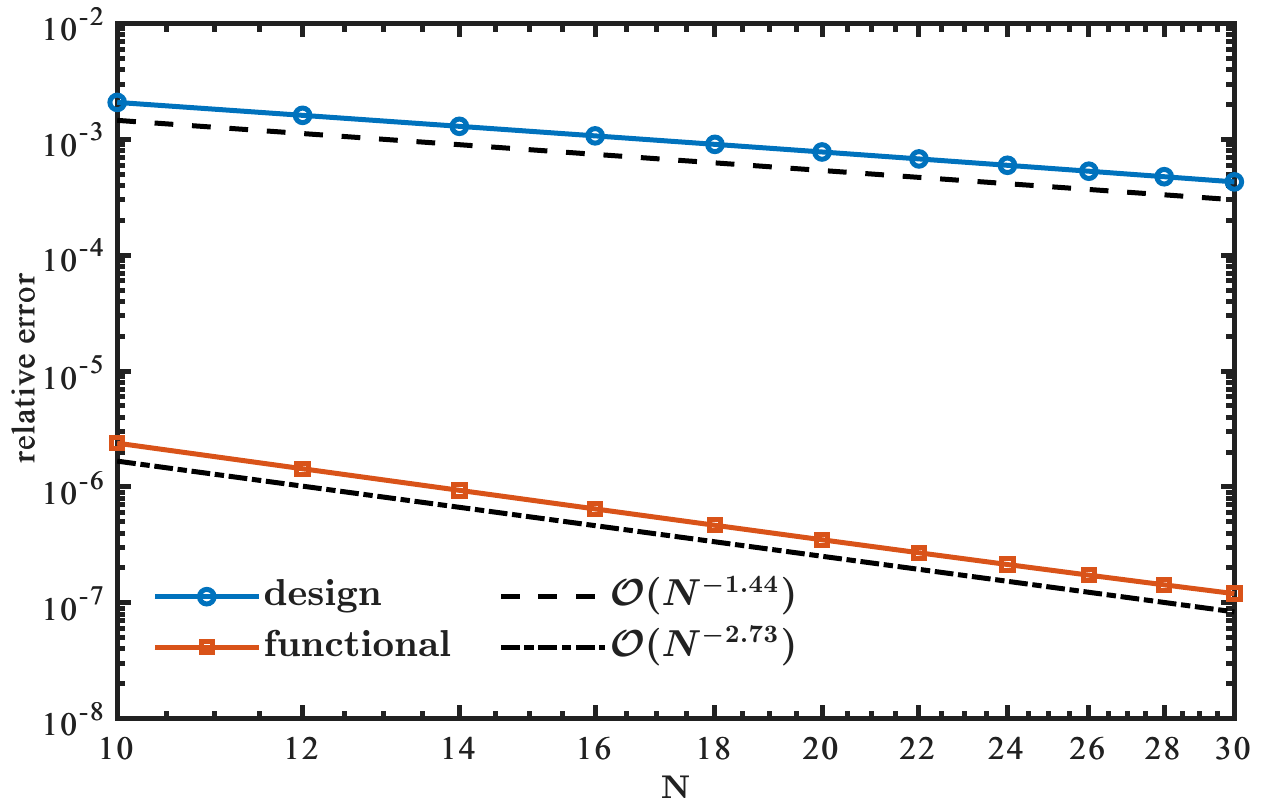}
	\caption{Convergence of the optimized design and optimal objective value under space--time refinement for ST--SE scheme.
    The design error decreases at an observed rate close to \(\mathcal{O}(N^{-3/2})\), while the optimal objective-value error decreases faster, close to \(\mathcal{O}(N^{-5/2})\)
    \(N := N_t = N_x\).
     }
	\label{fig:Func}
\end{figure}

Figure~\ref{fig:Func} compares the convergence of the discrete optimal design and the corresponding optimal objective value under space--time refinement. The optimal design error decreases at an observed rate of approximately \(3/2\), while the error in the optimal objective value decays faster, at an observed rate of approximately \(5/2\). This faster objective convergence is consistent with stationarity of the reduced optimization problem near the optimum. Since the primal solution is eliminated through \(\calR(\u,\vrho)=\zero\), the optimization is governed by the reduced objective
\(\Jac(\vrho).\)
Hence, formally,
\[\Jac(\vrho_N)-\Jac(\vrho_\star) = \mathcal{O}\!\left(\|\vrho_N-\vrho_\star\|^2\right). \]
Together, these numerical  results support the dual-consistent construction and indicate improved accuracy for scalar functional quantities, while a complete proof of the functional convergence behavior in the space--time topology optimization setting is left for future work.

\subsection{Computational time and accuracy in solving topology optimization problem}\label{sec:num:time}
\begin{figure}
	\centering
    \includegraphics[width=0.53\linewidth]{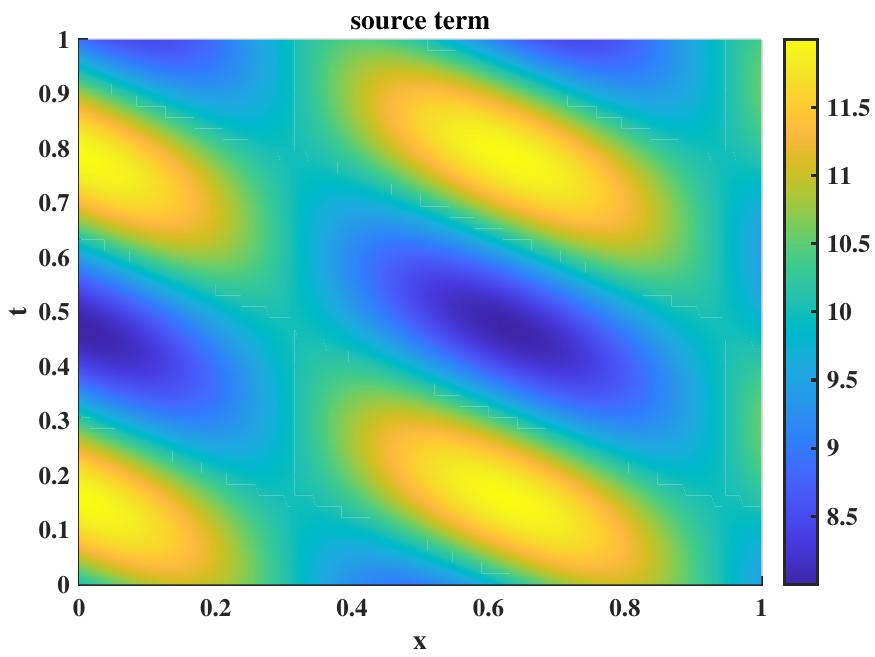}
    \caption{The source term $f(x,t)$ used in the numerical experiments.}
    \label{fig:source}
\end{figure}

We consider a model topology optimization problem subject to a 1D transient heat conduction problem on
\([0,1]\times[0,1]\), with a zero Neumann boundary condition at \(x=0\),
a zero Dirichlet boundary condition at \(x=1\), and zero initial condition:
\begin{align*}
\begin{aligned}
& u _t - (\kappa(x)  u _x)_x = 10+\sin\!\big(10(x+t)\big)+\sin(10\,t), &\qquad& (x,t)\in (0,1)\times(0,1),\\
& u_x (0,t) = 0,\quad  u (1,t)= 0,&\qquad& t\in[0,1 ],\\
& u (x,0) = 0,&\qquad& x\in[0,1].
\end{aligned}
\end{align*}
Figure \ref{fig:source} shows the source term $f(x,t)$ used in the numerical experiments. Both  FE and SE discretization employs $K=50$ spatial elements in space and a single global spectral element over the full time interval.
The FE discretization is based on linear shape functions, while the spectral element discretization uses 5 nodes per element (degree-4 polynomial basis), this number of nodes in the spatial direction is chosen to improve accuracy while keeping the cost reasonable.  Also in the construction of the ST–SE scheme, the homogeneous Neumann boundary condition is imposed weakly through an appropriate SAT term and all SAT parameters are chosen to satisfy the corresponding stability conditions. Further details are omitted for brevity.

In all runs, the volume constraint~\eqref{eq00} is enforced with a prescribed material budget \(V^{*}=0.5\), and we set \(p=3\) in~\eqref{eq0} as it is common in topology optimization, with $\kappa_{\min}=0.05$ and $\kappa_{\max}=1$. Although setting \(p=3\) penalizes intermediate densities in the material distribution, the density-based formulation may still yield non-binary designs, especially in one spatial dimension.

Topology optimization is performed with MMA, starting from the same initial design for all schemes at all time resolution levels. MMA iterations stop when the maximum design update per step falls below $10^{-4}$. Let \(\vrho^{(i)}\) and \(\Jac^{(i)}\) denote the optimized design and objective obtained at temporal resolution level \(i\) ($N^e_t$ for FE and $N^c_t$ for SE), we record the design change, 
\[\|\Delta\vrho\|:=\|\vrho^{(i)}-\vrho^{(i-1)}\|_\infty,\]
and the relative change of the objective,
\[\Jac_{\mbox{rel}}:=\,\frac{|\Jac^{(i)}-\Jac^{(i-1)}|}{\max(|\Jac^{(i)
|},\,|\Jac^{(i-1)}|,\,\varepsilon)},\]
 with \(\varepsilon=10^{-12}\).
We use two specified tolerances \(\tau_{\rho}\) for the design change and \(\tau_{\Jac}\) for the relative objective change.

Table~\ref{tab:timing} reports the wall-clock time and \(\|\Delta\vrho\|\) for each setting. Convergence is declared when $\|\Delta\vrho\|$ is below \(\tau_{\rho}=10^{-4}\) for two consecutive sweeps. The smallest design change occurs at \(N^e_t=8192\) for BE--FE--AAO, and at \(N^c_t=15\) for ST--SE. In a fair comparison between ST--SE and BE--FE--AAO, we deduce that ST--SE attains the desired tolerance with the fewest DoFs in time (accuracy-cost in terms of memory) and the shortest wall-clock time (time-to-solution).   Moreover, the spatial design resolution is determined by the number of elements $K$, whereas the accuracy of the state and adjoint solutions can be improved independently by increasing $N_x$ and $N_t$, without introducing additional design variables.

\begin{table}[htbp]
\centering
\small
\begin{tabular}{| rcc | rcc|}
\hline
\multicolumn{3}{|c|}{\textbf{BE--FE--AAO}} & 
\multicolumn{3}{c|}{\textbf{ST--SE}} \\
\hline
$N^e_t$ & time [s] & $\|\Delta\vrho\|$ & 
$N^c_t$ & time [s] & $\|\Delta\vrho\|$ \\
\hline
8     & 0.1007 & 0.050151  & 11 & 2.2695 & 0.000338 \\
16    & 0.1237 & 0.338284  & 13 & 2.4287 & 0.000025 \\
32    & 0.1949 & 0.011547  & 15 & 2.4863 & 0.000001 \\
64    & 0.3699 & 0.006566  &    &        &          \\
128   & 0.6206 & 0.003575  &    &        &          \\
256   & 1.4553 & 0.001871  &    &        &          \\
512   & 3.0582 & 0.000922  &    &        &          \\
1024  & 6.8655 & 0.000482  &    &        &          \\
2048  & 13.1422 & 0.000243 &    &        &          \\
4096  & 30.0512 & 0.000122 &    &        &          \\
8192  & 58.2616 & 0.000061 &    &        &          \\
\hline
\end{tabular}
\caption{Timing and design change comparison among the two all-at-once forward solvers. Left:  BE--FE--AAO (all-at-once version of BE--FE with varying $N_t^e$); right: ST--SE (space--time spectral elements with one time element, varying the number of temporal collocation nodes $N_t^c$). For the finest discretizations shown, the forward linear system sizes are $417843$ for BE--FE--AAO and $4000$ for ST--SE.}
\label{tab:timing}
\end{table}

\begin{figure}
	\centering
    \includegraphics[width=0.45\linewidth]{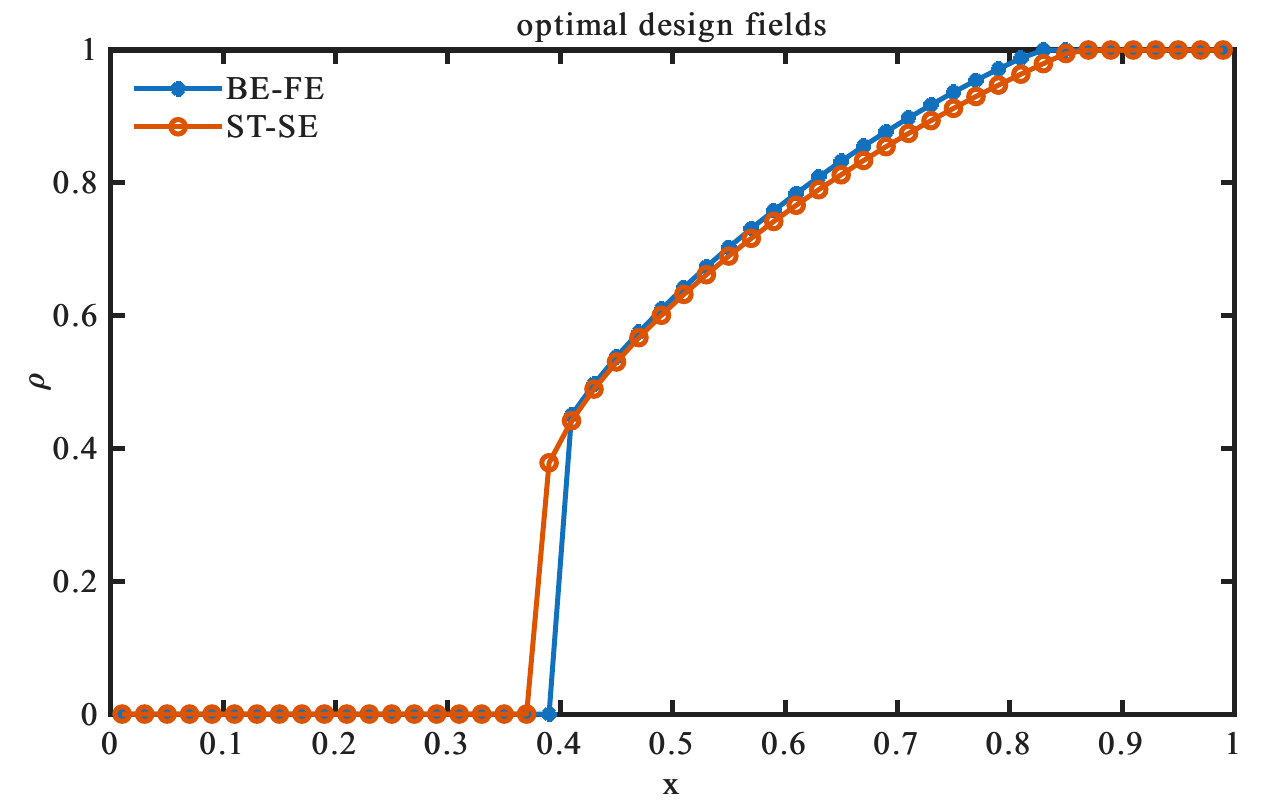}
    \includegraphics[width=0.45\linewidth]{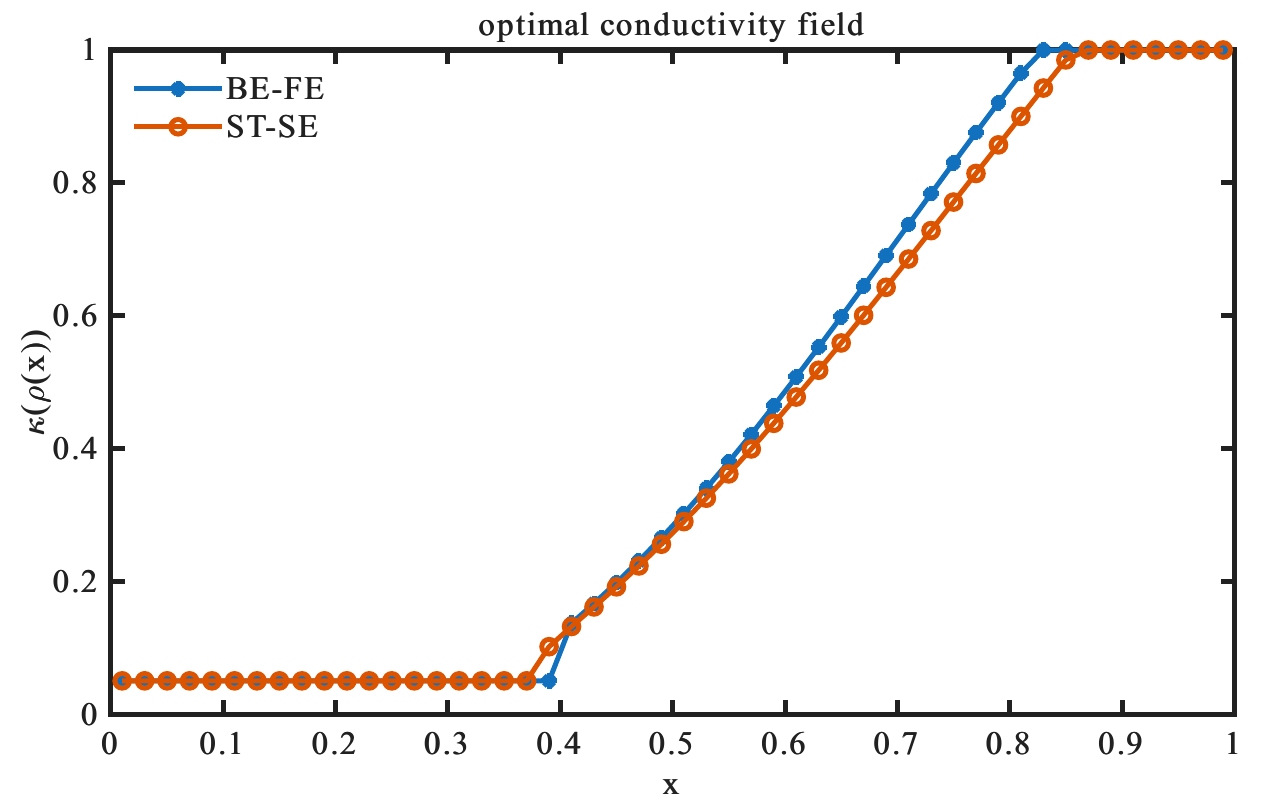}
    \caption{Comparison of optimal designs and diffusivity fields from the BE--FE--AAO and ST--SE approaches. }
    \label{fig:optDes}
\end{figure}
Figure~\ref{fig:optDes} shows the design field $\rho$ and thermal diffusivity $\kappa(\rho(x))$ at the finest resolution. The  transition region reflects the continuous density formulation, where intermediate densities may remain even when penalized material interpolation is used. Because the left boundary is insulated (zero Neumann boundary condition) while the right boundary is a cold sink (zero-temperature Dirichlet boundary condition), the optimal layout concentrates a single high-conductivity block next to the right boundary to channel heat out, with an almost linear tapering of the conductivity towards the low conductivity region furthest from the cold sink.

\begin{figure}
	\centering
    \includegraphics[width=0.45\linewidth]{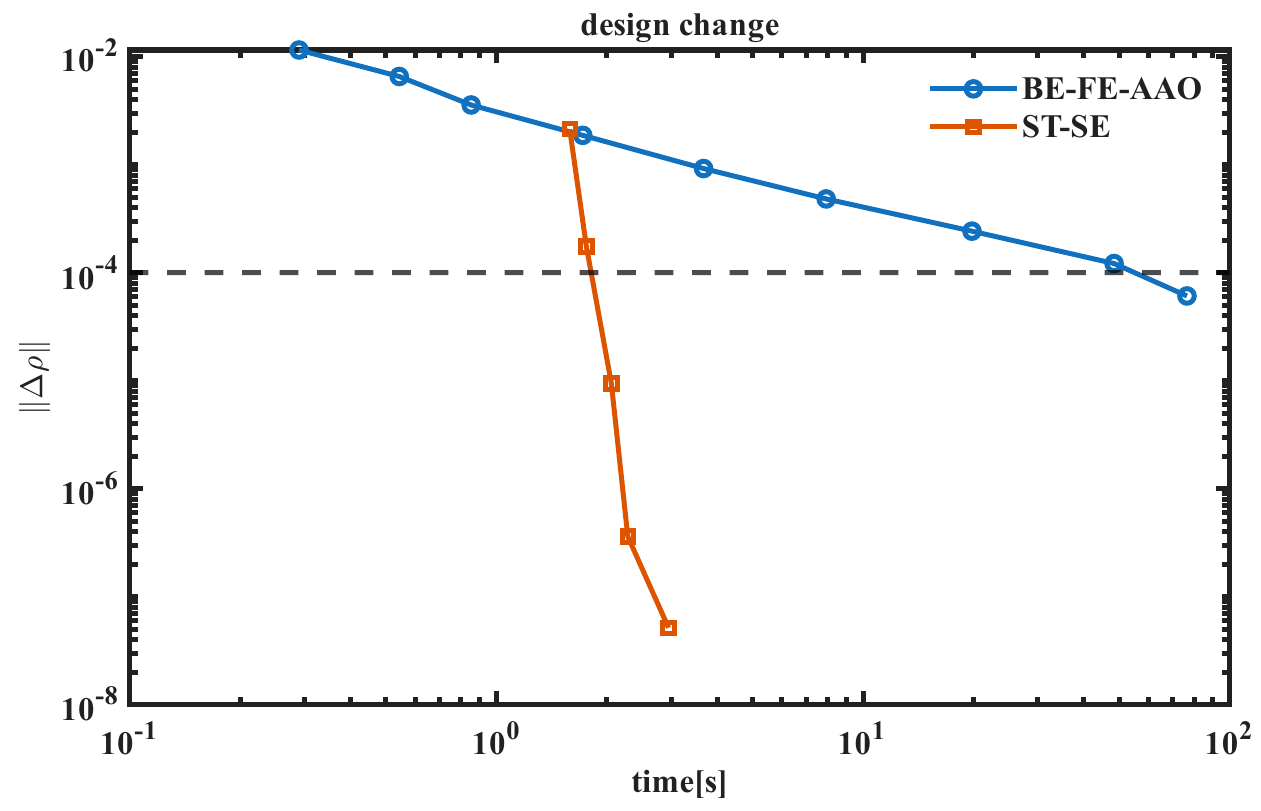}
    \includegraphics[width=0.45\linewidth]{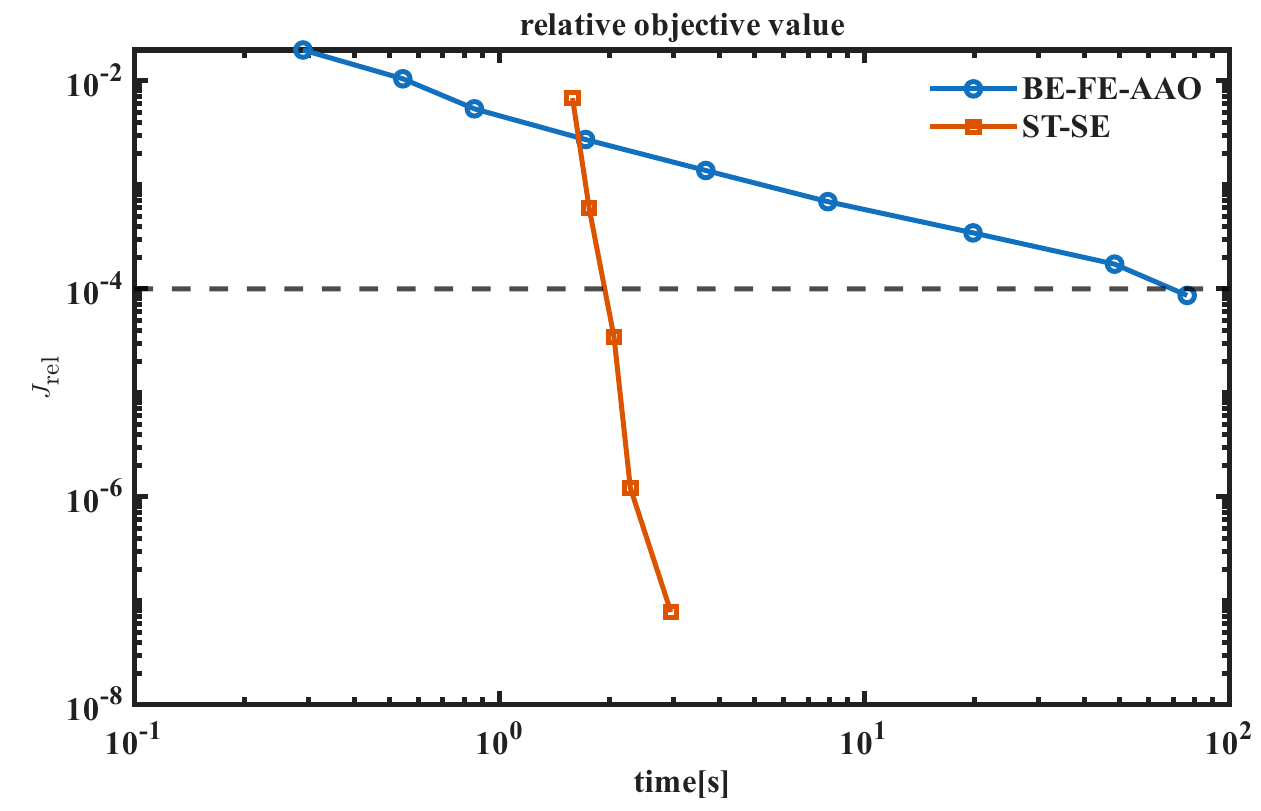}
\caption{Comparison of time from the two forward solvers 
(BE--FE--AAO and ST--SE) to reach the specified tolerance in the design change (left) and the relative change in objective (right).} 
    \label{fig:timeDes}
\end{figure}

Figure~\ref{fig:timeDes} shows the time required by the two forward solvers ST--SE and BE--FE-AAO, to reach the prescribed stopping tolerances for the design change. 
Comparing the two all-at-once methods, it is seen that ST--SE provides significantly higher accuracy at a lower computational time. In general all-at-once solvers associate  much larger computational systems, which requires parallelization to make practical sense and they become more attractive when coupled with space-time parallelization.

\begin{figure}[htbp]
	\centering
    \includegraphics[width=0.45\linewidth]{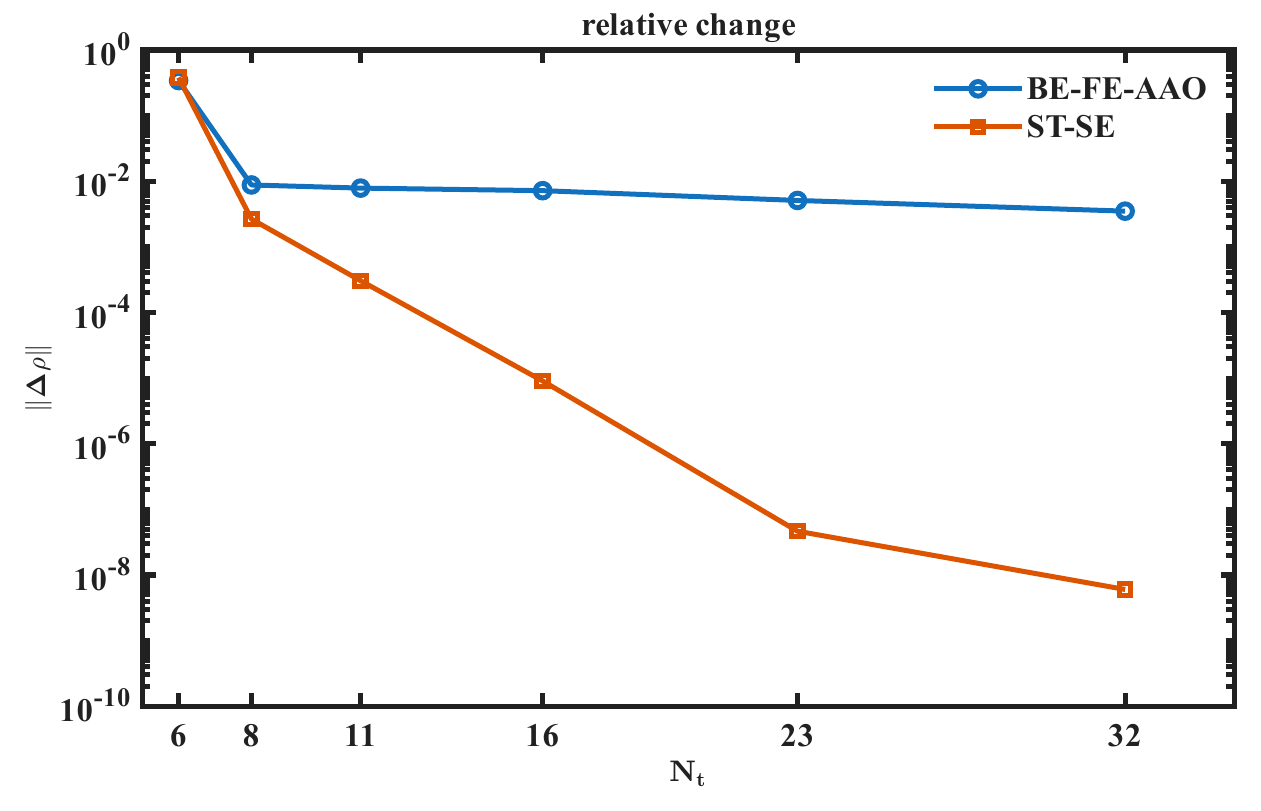}
    \includegraphics[width=0.45\linewidth]{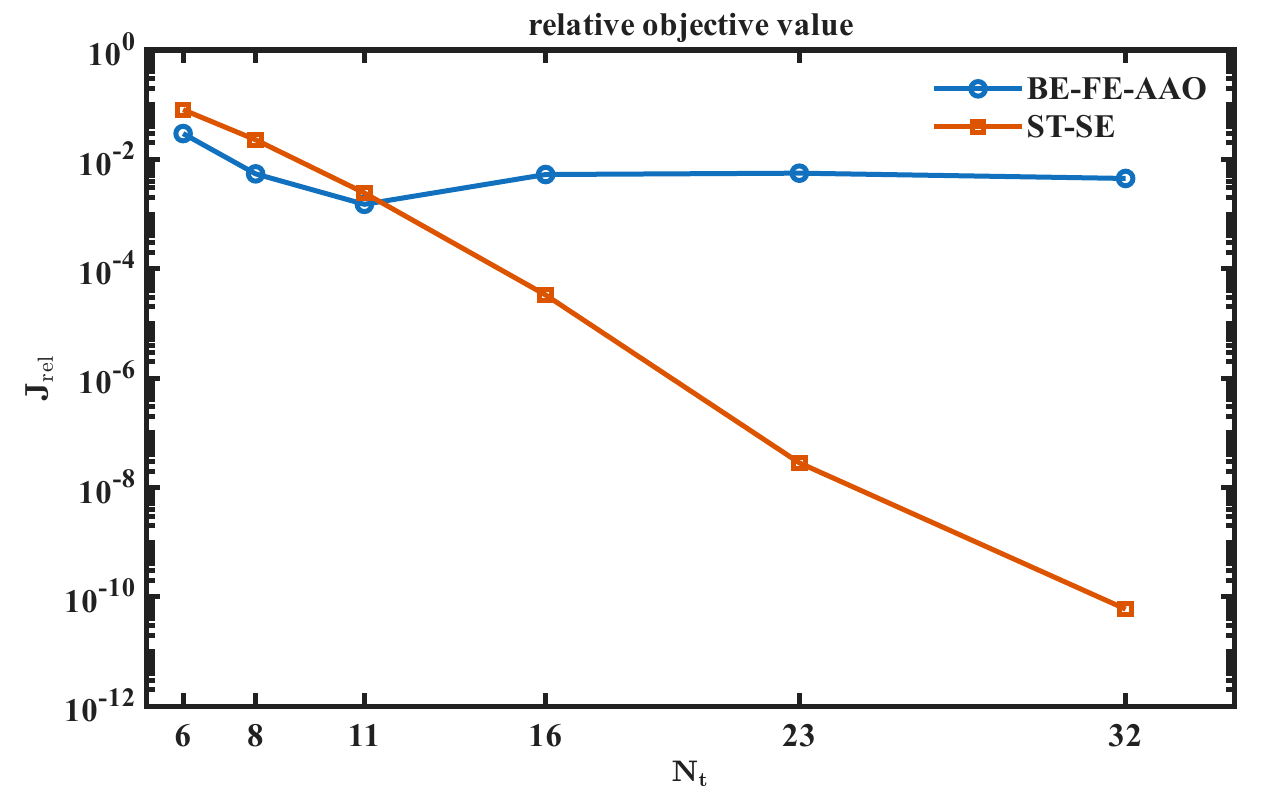}
    \caption{Comparison of design change (left)  and relative change of the objective (right) for two forward solvers BE--FE--AAO and ST--SE.}
    \label{fig:JreDes}
\end{figure}
Figure~\ref{fig:JreDes} demonstrates the evolution of the design change and the relative objective change over optimization iterations for BE--FE--AAO and ST--SE at the same number of nodes in temporal direction. The results demonstrate that the higher-order, dual-consistent space-time scheme offers a significantly more favorable accuracy–cost trade-off. It is clearly seen that the BE-FE-AAO scheme is converging at an exceptionally slow rate when compared to the proposed ST-SE scheme.

\section{Concluding remarks}\label{sec:sum}
We present a high-order space–time spectral element SBP–SAT framework for topology optimization of transient heat conduction. Interface and boundary conditions are imposed weakly via SATs, ensuring stability on heterogeneous subdomains. Using a discrete, dual-consistent space–time adjoint scheme, we compute design sensitivities over the full time horizon. Under standard smoothness assumptions, the resulting functional estimates exhibit spectral convergence behavior. Numerical experiments produce optimal designs consistent with independently computed reference solutions. Compared with low-order backward Euler finite-element all-at-once solvers, our space--time scheme achieves superior time-to-solution and accuracy-per-cost. In particular, compared with the BE--FE--AAO formulation, the proposed high-order space–time discretization attains higher resolution with fewer degrees of freedom, reducing computational time and memory while remaining stable. Large-scale performance further depends on effective linear solver technology. Some promising directions for future work include designing robust preconditioners for the coupled space–time linear systems and developing parallel-in-time strategies compatible with the monolithic space-time formulation.

\bibliographystyle{elsarticle-num} 
\bibliography{references}

@book {Lanczos1961,
    AUTHOR = {Lanczos, Cornelius},
     TITLE = {Linear differential operators},
 PUBLISHER = {D. Van Nostrand Co., Ltd., London-Toronto-New York-Princeton,
              N.J.},
      YEAR = {1961},
     PAGES = {xvi+564},
   MRCLASS = {34.00 (47.60)},
  MRNUMBER = {129153},
MRREVIEWER = {P.\ Roman},
URL={https://epubs.siam.org/doi/book/10.1137/1.9781611971187},
}

@article {Lindstrom2010,
	AUTHOR = {Lindstr\"{o}m, Jens and Nordstr\"{o}m, Jan},
	TITLE = {A stable and high-order accurate conjugate heat transfer
	problem},
	JOURNAL = {J. Comput. Phys.},
	FJOURNAL = {Journal of Computational Physics},
	VOLUME = {229},
	YEAR = {2010},
	NUMBER = {14},
	PAGES = {5440--5456},
	ISSN = {0021-9991,1090-2716},
	MRCLASS = {80A20 (65M06)},
	MRNUMBER = {2646148},
    DOI = {10.1016/j.jcp.2010.04.010},
}

@article {Carpenter2010,
	AUTHOR = {Carpenter, Mark H. and Nordstr\"{o}m, Jan and Gottlieb, David},
	TITLE = {Revisiting and extending interface penalties for multi-domain
	summation-by-parts operators},
	JOURNAL = {J. Sci. Comput.},
	FJOURNAL = {Journal of Scientific Computing},
	VOLUME = {45},
	YEAR = {2010},
	NUMBER = {1-3},
	PAGES = {118--150},
	ISSN = {0885-7474,1573-7691},
	MRCLASS = {65M60 (65M06 65M12)},
	MRNUMBER = {2679793},
	MRREVIEWER = {Istv\'{a}n\ Farag\'{o}},	
    DOI = {10.1007/s10915-009-9301-5},
}

@article {DelRey2018,
	AUTHOR = {Del Rey Fern\'andez, David C. and Hicken, Jason E. and Zingg,
	David W.},
	TITLE = {Simultaneous approximation terms for multi-dimensional
	summation-by-parts operators},
	JOURNAL = {J. Sci. Comput.},
	FJOURNAL = {Journal of Scientific Computing},
	VOLUME = {75},
	YEAR = {2018},
	NUMBER = {1},
	PAGES = {83--110},
	ISSN = {0885-7474,1573-7691},
	MRCLASS = {65M06 (65M12)},
	MRNUMBER = {3770313},
    DOI = {10.1007/s10915-017-0523-7},
}

@article {Carpenter1999,
	AUTHOR = {Carpenter, Mark H. and Nordstr\"om, Jan and Gottlieb, David},
	TITLE = {A stable and conservative interface treatment of arbitrary
	spatial accuracy},
	JOURNAL = {J. Comput. Phys.},
	FJOURNAL = {Journal of Computational Physics},
	VOLUME = {148},
	YEAR = {1999},
	NUMBER = {2},
	PAGES = {341--365},
	ISSN = {0021-9991,1090-2716},
	MRCLASS = {76M20},
	MRNUMBER = {1669703},
    DOI = {10.1006/jcph.1998.6114},
}

@article {Svard2014,
	AUTHOR = {Sv\"{a}rd, Magnus and Nordstr\"{o}m, Jan},
	TITLE = {Review of summation-by-parts schemes for
	initial-boundary-value problems},
	JOURNAL = {J. Comput. Phys.},
	FJOURNAL = {Journal of Computational Physics},
	VOLUME = {268},
	YEAR = {2014},
	PAGES = {17--38},
	ISSN = {0021-9991,1090-2716},
	MRCLASS = {65M06 (65M12)},
	MRNUMBER = {3192433},
    DOI = {10.1016/j.jcp.2014.02.031},
}

@article {DelRey2014,
	AUTHOR = {Del Rey Fern\'{a}ndez, David C. and Hicken, Jason E. and
	Zingg, David W.},
	TITLE = {Review of summation-by-parts operators with simultaneous approximation terms for the numerical solution of partial differential equations},
	JOURNAL = {Comput. \& Fluids},
	FJOURNAL = {Computers \& Fluids. An International Journal},
	VOLUME = {95},
	YEAR = {2014},
	PAGES = {171--196},
	ISSN = {0045-7930,1879-0747},
	MRCLASS = {65M06},
	MRNUMBER = {3189003},
	MRREVIEWER = {Sergio\ Amat},
    DOI = {10.1016/j.compfluid.2014.02.016},
}

@article {Mattsson2003,
    AUTHOR = {Mattsson, Ken},
     TITLE = {Boundary procedures for summation-by-parts operators},
   JOURNAL = {J. Sci. Comput.},
  FJOURNAL = {Journal of Scientific Computing},
    VOLUME = {18},
      YEAR = {2003},
    NUMBER = {1},
     PAGES = {133--153},
      ISSN = {0885-7474,1573-7691},
   MRCLASS = {65M06 (76M20)},
  MRNUMBER = {1958938},
   DOI = {10.1023/A:1020342429644},
}

@article {Nordstrom2013,
	AUTHOR = {Nordstr\"{o}m, Jan and Berg, Jens},
	TITLE = {Conjugate heat transfer for the unsteady compressible
	{N}avier-{S}tokes equations using a multi-block coupling},
	JOURNAL = {Comput. \& Fluids},
	FJOURNAL = {Computers \& Fluids. An International Journal},
	VOLUME = {72},
	YEAR = {2013},
	PAGES = {20--29},
	ISSN = {0045-7930,1879-0747},
	MRCLASS = {76N10 (65M06 65M12 76M20 80A20)},
	MRNUMBER = {3035902},
    DOI = {10.1016/j.compfluid.2012.11.018},
}

@article {Lundquist2018,
	AUTHOR = {Lundquist, Tomas and Malan, Arnaud and Nordstr\"{o}m, Jan},
	TITLE = {A hybrid framework for coupling arbitrary summation-by-parts
	schemes on general meshes},
	JOURNAL = {J. Comput. Phys.},
	FJOURNAL = {Journal of Computational Physics},
	VOLUME = {362},
	YEAR = {2018},
	PAGES = {49--68},
	ISSN = {0021-9991,1090-2716},
	MRCLASS = {65M08 (65M12 65M50)},
	MRNUMBER = {3774923},
	MRREVIEWER = {Ali\ R.\ Soheili},
	DOI = {10.1016/j.jcp.2018.02.018},
}

@article{Gassner2013,
	author = {Gassner, Gregor J.},
title = {A Skew-Symmetric Discontinuous Galerkin Spectral Element Discretization and Its Relation to SBP-SAT Finite Difference Methods},
journal = {SIAM Journal on Scientific Computing},
volume = {35},
number = {3},
pages = {A1233-A1253},
year = {2013},
doi = {10.1137/120890144},
}

@article{Hicken2011,
author = {Hicken, Jason E. and Zingg, David W.},
title = {Superconvergent Functional Estimates from Summation-By-Parts Finite-Difference Discretizations},
journal = {SIAM Journal on Scientific Computing},
volume = {33},
number = {2},
pages = {893-922},
year = {2011},
doi = {10.1137/100790987},
}

@article {Berg2012,
    AUTHOR = {Berg, Jens and Nordstr\"om, Jan},
     TITLE = {Superconvergent functional output for time-dependent problems
              using finite differences on summation-by-parts form},
   JOURNAL = {J. Comput. Phys.},
  FJOURNAL = {Journal of Computational Physics},
    VOLUME = {231},
      YEAR = {2012},
    NUMBER = {20},
     PAGES = {6846--6860},
      ISSN = {0021-9991,1090-2716},
   MRCLASS = {65M06 (65M12)},
  MRNUMBER = {2965104},
MRREVIEWER = {Kazufumi\ Ozawa},
       DOI = {10.1016/j.jcp.2012.06.032},
}

@article {Berg2013,
    AUTHOR = {Berg, Jens and Nordstr\"om, Jan},
     TITLE = {On the impact of boundary conditions on dual consistent finite
              difference discretizations},
   JOURNAL = {J. Comput. Phys.},
  FJOURNAL = {Journal of Computational Physics},
    VOLUME = {236},
      YEAR = {2013},
     PAGES = {41--55},
      ISSN = {0021-9991,1090-2716},
   MRCLASS = {65M06 (65M12)},
  MRNUMBER = {3020043},
       DOI = {10.1016/j.jcp.2012.11.019},
}

@article {Nikkar2019,
    AUTHOR = {Nikkar, Samira and Nordstr\"om, Jan},
     TITLE = {A dual consistent summation-by-parts formulation for the
              linearized incompressible {N}avier-{S}tokes equations posed on
              deforming domains},
   JOURNAL = {J. Comput. Phys.},
  FJOURNAL = {Journal of Computational Physics},
    VOLUME = {376},
      YEAR = {2019},
     PAGES = {322--338},
      ISSN = {0021-9991,1090-2716},
   MRCLASS = {65M06 (65M12 76D05 76M20)},
  MRNUMBER = {3875524},
       DOI = {10.1016/j.jcp.2018.09.006},
}

@article {Worku2022,
    AUTHOR = {Worku, Zelalem Arega and Zingg, David W.},
     TITLE = {Stability and functional superconvergence of narrow-stencil
              second-derivative generalized summation-by-parts
              discretizations},
   JOURNAL = {J. Sci. Comput.},
  FJOURNAL = {Journal of Scientific Computing},
    VOLUME = {90},
      YEAR = {2022},
    NUMBER = {1},
     PAGES = {Paper No. 42, 32},
      ISSN = {0885-7474,1573-7691},
   MRCLASS = {65M06 (65M12)},
  MRNUMBER = {4350188},
MRREVIEWER = {Swarn\ Singh},
       DOI = {10.1007/s10915-021-01707-5},
}

@article {Nordstrom2013time,
    AUTHOR = {Nordstr\"om, Jan and Lundquist, Tomas},
     TITLE = {Summation-by-parts in time},
   JOURNAL = {J. Comput. Phys.},
  FJOURNAL = {Journal of Computational Physics},
    VOLUME = {251},
      YEAR = {2013},
     PAGES = {487--499},
      ISSN = {0021-9991,1090-2716},
   MRCLASS = {65M20 (39A12 65B10 65L05)},
  MRNUMBER = {3094932},
MRREVIEWER = {Arnak\ Poghosyan},
       DOI = {10.1016/j.jcp.2013.05.042},
}

@article {Nordstrom2016,
    AUTHOR = {Nordstr\"om, Jan and Lundquist, Tomas},
     TITLE = {Summation-by-parts in time: the second derivative},
   JOURNAL = {SIAM J. Sci. Comput.},
  FJOURNAL = {SIAM Journal on Scientific Computing},
    VOLUME = {38},
      YEAR = {2016},
    NUMBER = {3},
     PAGES = {A1561--A1586},
      ISSN = {1064-8275,1095-7197},
   MRCLASS = {65M06 (65M12)},
  MRNUMBER = {3505309},
MRREVIEWER = {Pratibhamoy\ Das},
       DOI = {10.1137/15M103861X},
}

@article {Lucas2019,
    AUTHOR = {Friedrich, Lucas and Schn\"ucke, Gero and Winters, Andrew R.
              and Del Rey Fern\'andez, David C. and Gassner, Gregor J. and
              Carpenter, Mark H.},
     TITLE = {Entropy stable space-time discontinuous {G}alerkin schemes
              with summation-by-parts property for hyperbolic conservation
              laws},
   JOURNAL = {J. Sci. Comput.},
  FJOURNAL = {Journal of Scientific Computing},
    VOLUME = {80},
      YEAR = {2019},
    NUMBER = {1},
     PAGES = {175--222},
      ISSN = {0885-7474,1573-7691},
   MRCLASS = {65M60 (35L65 65M12)},
  MRNUMBER = {3954440},
MRREVIEWER = {B\"ulent\ Karas\"ozen},
       DOI = {10.1007/s10915-019-00933-2},
}

@article{Appel2024,
	author = {Appel, Magnus and Alexandersen, Joe},
title = {One-Shot Parareal Approach for Topology Optimization of Transient Heat Flow},
journal = {SIAM Journal on Scientific Computing},
volume = {47},
number = {6},
pages = {B1450-B1480},
year = {2025},
doi = {10.1137/24M1696603},
}

@article{Appel2025,
	title={Space-Time Multigrid Methods Suitable for Topology Optimisation of Transient Heat Conduction},
	author={Appel, Magnus and Alexandersen, Joe},
	journal={arXiv preprint arXiv:2505.10168},
	year={2025},
	url={https://arxiv.org/abs/2505.10168},
}

@Article{Alex2020,
AUTHOR = {Alexandersen, Joe and Andreasen, Casper Schousboe},
TITLE = {A Review of Topology Optimisation for Fluid-Based Problems},
JOURNAL = {Fluids},
VOLUME = {5},
YEAR = {2020},
NUMBER = {1},
ARTICLE-NUMBER = {29},
ISSN = {2311-5521},
DOI = {10.3390/fluids5010029}
}

@article{Alex2025,
	title={Large-Scale Topology Optimisation of Time-dependent Thermal Conduction Using Space-Time Finite Elements and a Parallel Space-Time Multigrid Preconditioner},
	author={Alexandersen, Joe and Appel, Magnus},
	journal={arXiv preprint arXiv:2508.09589}, 	
	year={2025},
	url={https://arxiv.org/abs/2508.09589},
}

@Article{Zeng2020,
  author  = {Zeng, Tao and Wang, Hu and Yang, Mengzhu and Alexandersen, Joe},
  journal = {International Journal of Heat and Mass Transfer},
  title   = {Topology optimization of heat sinks for instantaneous chip cooling using a transient pseudo-3D thermofluid model},
  year    = {2020},
  issn    = {00179310},
  volume  = {154},
  doi     = {10.1016/j.ijheatmasstransfer.2020.119681},
  type    = {Journal Article},
}

@article{Wu2019,
  title={Topology optimization for minimizing the maximum temperature of transient heat conduction structure},
  author={Wu, Shuhao and Zhang, Yongcun and Liu, Shutian},
  journal={Structural and Multidisciplinary Optimization},
  volume={60},
  pages={69--82},
  year={2019},
  publisher={Springer},
  doi={https://doi.org/10.1007/s00158-019-02196-9},
}

@article{Zhuang2013,
author = {Zhuang, Chungang and Xiong, Zhenhua and Ding, Han},
journal = {Numerical Heat Transfer. Part B, Fundamentals},
number = {3},
pages = {239-262},
title = {Topology Optimization of the Transient Heat Conduction Problem on a Triangular Mesh},
volume = {64},
year = {2013},
doi = {10.1080/10407790.2013.791785}
}

@article{Zhuang2014,
author = {Zhuang, Chungang and Xiong, Zhenhua},
journal = {Numerical Heat Transfer. Part B, Fundamentals},
number = {5},
pages = {445-471},
publisher = {Taylor & Francis Group},
title = {A Global Heat Compliance Measure Based Topology Optimization for the Transient Heat Conduction Problem},
volume = {65},
year = {2014},
doi = {10.1080/10407790.2013.873309}
}

@article {Nobis2022,
	AUTHOR = {Nobis, Harrison and Schlatter, Philipp and Wadbro, Eddie and
	Berggren, Martin and Henningson, Dan S.},
	TITLE = {Topology optimization of unsteady flows using the spectral
	element method},
	JOURNAL = {Comput. \& Fluids},
	FJOURNAL = {Computers \& Fluids. An International Journal},
	VOLUME = {239},
	YEAR = {2022},
	PAGES = {Paper No. 105387, 14},
	ISSN = {0045-7930,1879-0747},
	MRCLASS = {65M60 (76D05)},
	MRNUMBER = {4398362},
	DOI = {10.1016/j.compfluid.2022.105387},
}

@article {Nobis2023,
	AUTHOR = {Nobis, Harrison and Schlatter, Philipp and Wadbro, Eddie and
	Berggren, Martin and Henningson, Dan S.},
	TITLE = {Modal laminar-turbulent transition delay by means of topology
	optimization of superhydrophobic surfaces},
	JOURNAL = {Comput. Methods Appl. Mech. Engrg.},
	FJOURNAL = {Computer Methods in Applied Mechanics and Engineering},
	VOLUME = {403},
	YEAR = {2023},
	PAGES = {Paper No. 115721, 24},
	ISSN = {0045-7825,1879-2138},
	MRCLASS = {76F10},
	MRNUMBER = {4507502},
	DOI = {10.1016/j.cma.2022.115721},
}

@article{Saglietti2018,
	title={Topology optimization of heat sinks in a square differentially heated cavity},
	author={Saglietti, Clio and Schlatter, Philipp and Wadbro, Eddie and Berggren, Martin and Henningson, Dan S},
	journal={International Journal of Heat and Fluid Flow},
	volume={74},
	pages={36--52},
	year={2018},
	publisher={Elsevier},
    doi= {10.1016/j.ijheatfluidflow.2018.08.004},
}

@article{Yamaleev2010,
	title = {Local-in-time adjoint-based method for design optimization of unsteady flows},
	journal = {Journal of Computational Physics},
	volume = {229},
	number = {14},
	pages = {5394-5407},
	year = {2010},
	issn = {0021-9991},
	doi = {https://doi.org/10.1016/j.jcp.2010.03.045},
	author = {Nail K. Yamaleev and Boris Diskin and Eric J. Nielsen},
}

@article{Zahr2016,
	title = {An adjoint method for a high-order discretization of deforming domain conservation laws for optimization of flow problems},
	journal = {Journal of Computational Physics},
	volume = {326},
	pages = {516-543},
	year = {2016},
	issn = {0021-9991},
	doi = {https://doi.org/10.1016/j.jcp.2016.09.012},
	author = {M.J. Zahr and P.-O. Persson},
}

@article{Theulings2024,
	author = {Theulings, M. J. B. and Maas, R. and Noël, L. and van Keulen, F. and Langelaar, M.},
	title = {Reducing time and memory requirements in topology optimization of transient problems},
	journal = {International Journal for Numerical Methods in Engineering},
	volume = {125},
	number = {14},
	pages = {e7461},
	doi = {10.1002/nme.7461},
	year = {2024}
}

@article{Margetis2023,
	author = {Margetis, Andreas Stefanos I. and Papoutsis-Kiachagias, Evangelos M. and Giannakoglou, Kyriakos C.},
	title = {Reducing memory requirements of unsteady adjoint by synergistically using check-pointing and compression},
	journal = {International Journal for Numerical Methods in Fluids},
	volume = {95},
	number = {1},
	pages = {23-43},
	doi = {https://doi.org/10.1002/fld.5136},
	year = {2023}
}

@Article{Yaji2018,
	author={Yaji, Kentaro
	and Ogino, Masao
	and Chen, Cong
	and Fujita, Kikuo},
	title={Large-scale topology optimization incorporating local-in-time adjoint-based method for unsteady thermal-fluid problem},
	journal={Structural and Multidisciplinary Optimization},
	year={2018},
	month={Aug},
	day={01},
	volume={58},
	number={2},
	pages={817-822},
	doi={10.1007/s00158-018-1922-6},
}

@book{Sigmund2004,
  author    = {Bends{\o}e, Martin P. and Sigmund, Ole},
  title     = {Topology Optimization: theory, methods, and applications},
  edition   = {2},
  publisher = {Springer},
  address   = {Berlin Heidelberg},
  year      = {2004},
  pages     = {xiv+370},
  doi       = {10.1007/978-3-662-05086-6},
  mrclass   = {74P15 (49Q10 74P05 74Q05)},
  mrnumber  = {2008524},
  mrreviewer= {Michal Ko\v{c}vara}
}

@article{Bendsoe1988,
title = {Generating optimal topologies in structural design using a homogenization method},
journal = {Computer Methods in Applied Mechanics and Engineering},
volume = {71},
number = {2},
pages = {197-224},
year = {1988},
issn = {0045-7825},
doi = {https://doi.org/10.1016/0045-7825(88)90086-2},
author = {Martin Philip Bends{\o}e and Noboru Kikuchi},
}

@article {Sigmund2013,
    AUTHOR = {Sigmund, Ole and Maute, Kurt},
     TITLE = {Topology optimization approaches},
   JOURNAL = {Struct. Multidiscip. Optim.},
  FJOURNAL = {Structural and Multidisciplinary Optimization},
    VOLUME = {48},
      YEAR = {2013},
    NUMBER = {6},
     PAGES = {1031--1055},
      ISSN = {1615-147X,1615-1488},
   MRCLASS = {74P15 (49Q12)},
  MRNUMBER = {3138124},
       DOI = {10.1007/s00158-013-0978-6},
}

@article {Deaton2014,
    AUTHOR = {Deaton, Joshua D. and Grandhi, Ramana V.},
     TITLE = {A survey of structural and multidisciplinary continuum
              topology optimization: post 2000},
   JOURNAL = {Struct. Multidiscip. Optim.},
  FJOURNAL = {Structural and Multidisciplinary Optimization},
    VOLUME = {49},
      YEAR = {2014},
    NUMBER = {1},
     PAGES = {1--38},
      ISSN = {1615-147X,1615-1488},
   MRCLASS = {74P15},
  MRNUMBER = {3182450},
       DOI = {10.1007/s00158-013-0956-z},
}

@article{Dbouk2017,
title = {A review about the engineering design of optimal heat transfer systems using topology optimization},
journal = {Applied Thermal Engineering},
volume = {112},
pages = {841-854},
year = {2017},
issn = {1359-4311},
doi = {https://doi.org/10.1016/j.applthermaleng.2016.10.134},
author = {T. Dbouk},
}

@article {Svanberg1987,
    AUTHOR = {Svanberg, Krister},
     TITLE = {The method of moving asymptotes---a new method for structural
              optimization},
   JOURNAL = {Internat. J. Numer. Methods Engrg.},
  FJOURNAL = {International Journal for Numerical Methods in Engineering},
    VOLUME = {24},
      YEAR = {1987},
    NUMBER = {2},
     PAGES = {359--373},
      ISSN = {0029-5981,1097-0207},
   MRCLASS = {73K40},
  MRNUMBER = {875307},
       DOI = {10.1002/nme.1620240207},
}

@article {Niles2000,
    AUTHOR = {Pierce, Niles A. and Giles, Michael B.},
     TITLE = {Adjoint recovery of superconvergent functionals from {PDE}
              approximations},
   JOURNAL = {SIAM Rev.},
  FJOURNAL = {SIAM Review},
    VOLUME = {42},
      YEAR = {2000},
    NUMBER = {2},
     PAGES = {247--264},
      ISSN = {0036-1445,1095-7200},
   MRCLASS = {65N12 (76M12 76M20)},
  MRNUMBER = {1778357},
MRREVIEWER = {Dennis\ C.\ Jespersen},
       DOI = {10.1137/S0036144598349423},
}

@article {Cockburn2007,
    AUTHOR = {Cockburn, Bernardo and Ichikawa, Ryuhei},
     TITLE = {Adjoint recovery of superconvergent linear functionals from
              {G}alerkin approximations. {T}he one-dimensional case},
   JOURNAL = {J. Sci. Comput.},
  FJOURNAL = {Journal of Scientific Computing},
    VOLUME = {32},
      YEAR = {2007},
    NUMBER = {2},
     PAGES = {201--232},
      ISSN = {0885-7474,1573-7691},
   MRCLASS = {65N30 (65N15)},
  MRNUMBER = {2320570},
MRREVIEWER = {Emmanuel\ Creus\'e},
       DOI = {10.1007/s10915-007-9129-9},
}

@article {Giles2002,
    AUTHOR = {Giles, Michael B. and S\"uli, Endre},
     TITLE = {Adjoint methods for {PDE}s: a posteriori error analysis and
              postprocessing by duality},
   JOURNAL = {Acta Numer.},
  FJOURNAL = {Acta Numerica},
    VOLUME = {11},
      YEAR = {2002},
     PAGES = {145--236},
      ISSN = {0962-4929,1474-0508},
   MRCLASS = {65N15 (35G15 35J25 65N30)},
  MRNUMBER = {2009374},
MRREVIEWER = {Lucas\ J\'odar},
       DOI = {10.1017/S096249290200003X},
}

@article{Michaleris1994,
   author = {Michaleris, Panagiotis and Tortorelli, Daniel A. and Vidal, Creto A.},
   title = {Tangent operators and design sensitivity formulations for transient non-linear coupled problems with applications to elastoplasticity},
   journal = {International Journal for Numerical Methods in Engineering},
   volume = {37},
   number = {14},
   pages = {2471-2499},
   DOI = {10.1002/nme.1620371408},
   year = {1994},
   type = {Journal Article}
}

@Article{Dahl2008,
  author  = {Dahl, Jonas and Jensen, Jakob S. and Sigmund, Ole},
  journal = {Structural and Multidisciplinary Optimization},
  title   = {Topology optimization for transient wave propagation problems in one dimension},
  year    = {2008},
  issn    = {1615-1488},
  number  = {6},
  pages   = {585-595},
  volume  = {36},
  doi     = {10.1007/s00158-007-0192-5},
  type    = {Journal Article},
}

@Article{Giles2000,
author={Giles, Michael B.
and Pierce, Niles A.},
title={An Introduction to the Adjoint Approach to Design},
journal={Flow, Turbulence and Combustion},
year={2000},
month={Dec},
day={01},
volume={65},
number={3},
pages={393-415},
doi={10.1023/A:1011430410075},
}

@article {Hartmann2007,
    AUTHOR = {Hartmann, Ralf},
     TITLE = {Adjoint consistency analysis of discontinuous {G}alerkin
              discretizations},
   JOURNAL = {SIAM J. Numer. Anal.},
  FJOURNAL = {SIAM Journal on Numerical Analysis},
    VOLUME = {45},
      YEAR = {2007},
    NUMBER = {6},
     PAGES = {2671--2696},
      ISSN = {0036-1429,1095-7170},
   MRCLASS = {65N30},
  MRNUMBER = {2361907},
MRREVIEWER = {Koffi\ B.\ Fadimba},
       DOI = {10.1137/060665117},
}

@article {Arnold2001,
    AUTHOR = {Arnold, Douglas N. and Brezzi, Franco and Cockburn, Bernardo and Marini, L. Donatella},
     TITLE = {Unified analysis of discontinuous {G}alerkin methods for elliptic problems},
   JOURNAL = {SIAM J. Numer. Anal.},
  FJOURNAL = {SIAM Journal on Numerical Analysis},
    VOLUME = {39},
      YEAR = {2001},
    NUMBER = {5},
     PAGES = {1749--1779},
      ISSN = {0036-1429,1095-7170},
   MRCLASS = {65N30},
  MRNUMBER = {1885715},
       DOI = {10.1137/S0036142901384162},
}

@article{Nataj2026CHT,
  author       = {Nataj, Sarah and {Del Rey Fern{\'a}ndez}, David C. and Brown, David and Jaiman, Rajeev},
  title        = {A finite-difference summation-by-parts, conditionally stable partitioned algorithm for conjugate heat transfer problems},
  journal      = {arXiv:2602.17843},
  year         = {2026},
  url          = {https://arxiv.org/abs/2602.17843}
}

@article {Long2018,
    AUTHOR = {Long, Kai and Wang, Xuan and Gu, Xianguang},
     TITLE = {Multi-material topology optimization for the transient heat
              conduction problem using a sequential quadratic programming
              algorithm},
   JOURNAL = {Eng. Optim.},
  FJOURNAL = {Engineering Optimization},
    VOLUME = {50},
      YEAR = {2018},
    NUMBER = {12},
     PAGES = {2091--2107},
      ISSN = {0305-215X,1029-0273},
   MRCLASS = {74P15 (80A20 90C55)},
  MRNUMBER = {3863733},
       DOI = {10.1080/0305215X.2017.1417401},
}

@article {Li2023,
    AUTHOR = {Li, Shuai and Zhang, Yongcun and Liu, Shutian},
     TITLE = {Multi-material topology optimization method for transient
              thermal structure with constraint on permissible temperatures},
   JOURNAL = {Internat. J. Numer. Methods Engrg.},
  FJOURNAL = {International Journal for Numerical Methods in Engineering},
    VOLUME = {124},
      YEAR = {2023},
    NUMBER = {18},
     PAGES = {4022--4057},
      ISSN = {0029-5981,1097-0207},
   MRCLASS = {74P15 (74A15)},
  MRNUMBER = {4631853},
       DOI = {10.1002/nme.7299},
       URL = {https://doi-org.proxy1-bib.sdu.dk/10.1002/nme.7299},
}

\appendix

\section{A two--subdomain manufactured solution for the heat equation with constant heat source}\label{app:twoDomSol}
We consider the 1D heat equation with a constant source term $f\in\mathbb{R}$,
\begin{equation}\label{eq:pde}
 u _t = (\kappa(x)\, u_x)_{x} + f, \qquad (x,t)\in [0,1]\times[0,T],
\end{equation}
and a material interface at $x=\xi\in(0,1)$. The thermal diffusivity is piecewise constant,
\[
\kappa(x)=
\begin{cases}
\kappa_1, & 0\le x\le \xi,\\[2pt]
\kappa_2, & \xi< x\le 1,
\end{cases}
\qquad \kappa_1,\,\kappa_2>0.
\]
The Dirichlet boundary conditions are given as
\[
 u (0,t)=0,\qquad  u (1,t)= u _R,\quad t\in[0,T],
\]
where $ u _R$ is a given constant. The standard interface conditions at $x=\xi$ are given by
\[
 u (\xi^-,t)= u (\xi^+,t), \qquad
\kappa_1\, u _x(\xi^-,t)=\kappa_2\, u _x(\xi^+,t), \quad t\in[0,T].
\]
We seek a manufactured  solution $ u (x,t)$ to \eqref{eq:pde} as a combination of a steady part $ u _s(x)$ and a decaying transient part $ v (x,t)$:
\[
 u (x,t)= u _s(x)+ v (x,t).
\]
We choose $ u _s(x)$ to absorb the constant source term and enforce the Dirichlet boundary conditions, and $ v (x,t)$ to enforce time dependence while keeping the boundary values fixed. Then the initial condition is given by \[ u (x,0)= u _s(x)+ v (x,t): = u _0(x), \quad x\in[0,1].\]
We calculate the steady part of the solution in the following way. In each subdomain we solve $\kappa_i\, u _s''+f=0$ with interface and boundary conditions. A convenient choice is a piecewise quadratic:
\[
u_s(x)=
\begin{cases}
-\dfrac{f}{2\kappa_1}\,x^2 + A_1\,x, & 0\le x\le \xi,\\[1.2ex]
-\dfrac{f}{2\kappa_2}\,x^2 + A_2\,x + B_2, & \xi\le x\le 1.
\end{cases}
\]
Imposing boundary and interface conditions gives the coefficients explicitly
\begin{align}
A_1 &= \frac{\;\kappa_2\, u _R\;+\;\dfrac{f}{2}\Big(1+\xi^2\big(\dfrac{\kappa_2}{\kappa_1}-1\big)\Big)\;}
{\;\xi\,\kappa_2 + (1-\xi)\,\kappa_1\;}, \label{eq:A1}\\[4pt]
A_2 &= \frac{\kappa_1}{\kappa_2}\,A_1,\qquad
B_2 \;=\;  u _R + \frac{f}{2\kappa_2} - A_2. \label{eq:A2B2}
\end{align}
By construction, $\kappa_i\, u _s''+f=0$ identically in each subdomain and all boundary and interface conditions hold.

Next  we calculate the transient part of the solution. 
Set $ v (x,t)=e^{-\lambda t}\, w (x)$ with $\lambda>0$ to be determined. Then $ v _t=\kappa(x)  v _{xx}$ if and only if
\[
\kappa_i\, w _i'' + \lambda\, w _i = 0 \quad \text{in each subdomain},
\]
with homogeneous Dirichlet conditions $ w (0)= w (1)=0$.
One form that satisfies this is
\begin{align}\label{eq:we}
    \begin{aligned}
        w_1(x)&=\sin(\alpha_1 x),\qquad &0\le x\le \xi,\\
       w_2(x)&=r\,\sin\big(\alpha_2(1-x)\big),\qquad & \xi\le x\le 1,   
    \end{aligned}
\end{align}
where $\alpha_i=\sqrt{\lambda/\kappa_i}$ and $r$ is chosen for temperature continuity:
\[
r=\frac{\sin(\alpha_1\xi)}{\sin(\alpha_2(1-\xi))}.
\]
The flux continuity at $x=\xi$ reads 
$\kappa_1  w _1'(\xi)=\kappa_2  w _2'(\xi)$,
producing  the following condition for $\lambda$:
\begin{equation}\label{eq:eig}
\;\sqrt{\kappa_1}\,\cot\!\big(\alpha_1\xi\big)\;+\;\sqrt{\kappa_2}\,\cot\!\big(\alpha_2(1-\xi)\big)\;=\;0,\;
\qquad \alpha_i=\sqrt{\lambda/\kappa_i}.
\end{equation}
Any positive root $\lambda$ of \eqref{eq:eig} yields a valid decaying transient.

Combining steady solution $u_s$ and transient correction  gives the manufactured solution
\begin{equation}\label{eq:finalu}
\; u (x,t)= u _s(x)+e^{-\lambda t}\, w (x)\;
\end{equation}
with $ u _s$ from \eqref{eq:A1}–\eqref{eq:A2B2}, $ w $ given by \eqref{eq:we}, and $\lambda$ solving \eqref{eq:eig}.

To verify the solution, we plug in \eqref{eq:finalu} in right-hand side of \eqref{eq:pde} in each subdomain,
\[
 u _t-\kappa_i  u _{xx} = \big(-\lambda e^{-\lambda t} w \big) - \kappa_i\big( u _s'' + e^{-\lambda t} w ''\big)
= -\lambda e^{-\lambda t} w  - \kappa_i\Big(-\frac{f}{\kappa_i}-\frac{\lambda}{\kappa_i}e^{-\lambda t} w \Big)=f.
\]
Thus \eqref{eq:pde}, Dirichlet boundary conditions and the interface conditions at $x=\xi$ hold.

\section{Adjoint equation for model topology optimization problem}\label{App:AdjSol}
Let \(x\in[0,1]\) and \(t\in[0,T]\), and set \(\Omega=(0,1)\times(0,T)\) and \( u \in H^{1}(\Omega)\).
Assume  $f,\,  u _R$ and $ u _0(x)$ are given.
Suppose there is an interface at \(\xi\in(0,1)\).
Define the primal PDE as follows:
\begin{align*}
\begin{aligned}
&\L\;  u  :=  u _t - (\kappa(x)\,  u _x)_x = f, &\qquad&(x,t)\in\Omega,\\
& u (0,t) = 0,\quad u(1,t)=u_R,&\qquad&   t\in[0,T],\\
& u (x,0) =  u_0(x),&\qquad&   x\in[0,1],\\
& u (\xi^-,t) =  u (\xi^+,t),\qquad\kappa_1  u _x(\xi^-,t)=\kappa_2  u _x(\xi^+,t),
\end{aligned}
\end{align*}
with piecewise positive $\kappa(x)$,
\[
\kappa(x)=
\begin{cases}
k_1,& 0\le x\le \xi,\\
k_2,& \xi<x\le 1,
\end{cases}
\qquad k_1,\, k_2>0.
\]
Assume the objective functional is given by
\[J( u )=\int_0^T\!\!\int_0^1  u (x,t)^2\,dx\,dt,\] 
Introduce the adjoint variable $ v \in H^1(\Omega)$  and the adjoint operator $\L^{*}$ satisfying \( ( v ,\L u )= (\L^{*}  v ,  u )\) by Green's identity \cite{Lanczos1961}.  The Lagrangian is given by
\[
\mathcal L(u,v)
=
J(u)
-
\int_0^T\int_0^1
v\big[u_t-(\kappa u_x)_x-f\big]\,dx\,dt := J(u)- N(u) .
\]
For a functional $F( u )$, its variation in the direction $\phi$ is
\[
\delta F[ u ;\phi] \;:=\; \left.\frac{d}{d\varepsilon}\,F( u +\varepsilon\phi)\right|_{\varepsilon=0},
\]
called the G\^ateaux derivative. Here $\phi$ is any test function satisfying the constraints where $u$ is fixed:
\[
\phi(0,t)=\phi(1,t)=0\quad, \qquad \phi(x,0)=0.
\]
and 
\[
\phi(\xi^-,t)=\phi(\xi^+,t),\qquad
\kappa_1\phi_x(\xi^-,t)=\kappa_2\phi_x(\xi^+,t).
\]
In the following we calculate $\partial L$. First note that the objective variation is given by
\[
\delta J[ u ;\phi]
=\left.\frac{d}{d\varepsilon}\int_0^T\!\!\int_0^1(u+\varepsilon\phi)^2 dx\, dt\right|_{\varepsilon=0}
=\int_0^T\!\!\int_0^1 2\,u\,\phi\,dx\,dt.
\]
The PDE variation is
\[
\delta N[ u ;\phi]
=\int_0^T\!\!\int_0^1  v \,\phi_t\,dx\,dt
\;-\;\int_0^T\!\!\int_0^1  v \,(\kappa\phi_x)_x\,dx\,dt.
\]
Integrate by parts in time:
\[
\int_0^T\!\!\int_0^1 v\,\phi_t
=\Big[\int_0^1 v\,\phi\,dx\Big]_{t=0}^{t=T}
-\int_0^T\!\!\int_0^1 v_t\,\phi\,dx\,dt
= \int_0^1 v(x,T)\,\phi(x,T)\,dx-\int_0^T\!\!\int_0^1 v_t\,\phi.
\]
We used $\phi(\cdot,0)=0$. 
Integrate by parts in space, piecewise on $(0,\xi)$ and $(\xi,1)$ (where $\kappa$ is constant):
\[
-\int_0^1 v\,(\kappa\phi_x)_x
= -\big[\kappa v\,\phi_x\big]_{\text{bdry/intf}}
+\int_0^1 \kappa v_x\,\phi_x
= -\big[\kappa v\,\phi_x\big]_{\text{bdry/intf}}
+\big[\kappa v_x\,\phi\big]_{\text{bdry/intf}}
-\int_0^1 (\kappa v_x)_x\,\phi.
\]
Using boundary and interface conditions, the PDE variation contributes
\[
\delta N[u;\phi]
=-\int_0^T\!\!\int_0^1 \big[v_t+(\kappa v_x)_x\big]\phi\,dx\,dt
+ \displaystyle \int_0^1 v(x,T)\phi(x,T)\,dx.\]
Therefore
\[
\delta\mathcal L
= \int_0^T\!\!\int_0^1 \big(2\,u + v_t + (\kappa v_x)_x\big)\,\phi\,dx\,dt
\;-\;\int_0^1 v(x,T)\,\phi(x,T)\,dx.
\]
Since $\phi$ is arbitrary in the interior and at $t=T$, the condition $\delta\mathcal L=0$ yields 
\[2\,u + v_t + (\kappa v_x)_x=0,\qquad v(x,T)=0. \]
Then the adjoint equation  is given by 
\begin{align*}
\begin{aligned}
&\calL^{*} v := -\,v_t - (\kappa v_x)_x = 2\,u, &\qquad&(x,t)\in\Omega,\\
&v(0,t)=v(1,t)=0,&\qquad&   t\in[0,T],\\
&v(x,T)=0,&\qquad&   x\in[0,1],\\
&v(\xi^-,t)=v(\xi^+,t),\quad \kappa_1 v_x(\xi^-,t)=\kappa_2 v_x(\xi^+,t).
\end{aligned}
\end{align*}

\end{document}